\documentclass[12pt,oneside]{amsart}
\usepackage{tikz}
\usetikzlibrary{shapes.geometric}
\usepackage{geometry}
\usepackage[latin1]{inputenc}
\usepackage{amsmath}
\usepackage{amsfonts}
\usepackage{amssymb}
\usepackage{graphicx}
\usepackage{epstopdf}

\newtheorem{thm}{Theorem}[section]
\newtheorem{asm}{Assumption}[section]
\newtheorem{defn}{Definition}[section]

\newtheorem{theorem}{Theorem}[section]
\newtheorem{corollary}[theorem]{Corollary}

\newtheorem{remark}{Remark}[theorem]

\def\bc{\begin{theorem}}
\def\ec{\end{theorem}}
\def\bc{\begin{corollary}}
\def\ec{\end{corollary}}
\def\be{\begin{equation}}
\def\ee{\end{equation}}
\def\bast{\begin{eqnarray*} }
\def\east{\end{eqnarray*} }
\def\bea{\begin{eqnarray}}
\def\eea{\end{eqnarray}}
\def\basm{\begin{asm}}
\def\easm{\end{asm}}

\def\ss{\smallskip}
\def\ms{\medskip}
\def\bs{\bigskip}
\def\ni{\noindent}

\def\lb{\linebreak}

\def\det{\hbox{det }}
\def\R{\mathbb R}

\def\ts{\tilde{s}}
\def\tr{\tilde{r}}
\def\ty{\tilde{y}}
\def\ww{\tilde\omega}

\def\D{\mathbb D}
\def\F{\mathcal F}

\def\tsig{\tilde{\sigma}}

\def\la{\langle}
\def\ra{\rangle}

\newcommand{\pder}[2]{\frac {\partial {#1}}{\partial {#2}}}

\def\dim{\textrm{ dim\,}}
\def\Ker{\textrm { ker\,}}

\def\rank{\textrm{ rank\,}}
\def\dim{\textrm{ dim\,}}
\def\R{\mathbb R}

\title[Microlocal Analysis of Borehole Seismic  Data]{Microlocal Analysis of Borehole Seismic  Data}
\author[Raluca Felea, Romina Gaburro, Allan Greenleaf and  Clifford Nolan]
{Raluca Felea, Romina Gaburro, Allan Greenleaf and  Clifford Nolan}


\address{School of Mathematical Sciences, Rochester Institute of Technology,  Rochester, NY, 14623}
\email{rxfsma@rit.edu}
\address{Department of Mathematics and Statistics, Health Research Institute (HRI), \ University of Limerick, Limerick, V94 T9PX, Ireland.}
\email{Romina.Gaburro@ul.ie}
\address{Department of Mathematics, University of Rochester, Rochester, NY, 14627}
\email{allan@math.rochester.edu}
\address{Department of Mathematics and Statistics, Health Research Institute (HRI), \ University of Limerick, Limerick, V94 T9PX, Ireland.}
\email{Clifford.Nolan@ul.ie}

\date{4 October, 2021.}

\begin{document}
\maketitle

\begin{abstract} Borehole seismic data is obtained by receivers located in a well, 
with sources located on the 
surface or in another well. Using microlocal analysis, 
we study possible approximate reconstruction via linearized, filtered 
backprojection of an isotropic  sound speed
in the subsurface for three types of data sets. The sources may form a dense array on the surface, 
or be located along a line on the surface (walkaway geometry) or in another borehole (crosswell).
We show that for the dense array, reconstruction   is feasible, with no artifacts in the absence of caustics in the 
background ray geometry, and mild artifacts in the presence of fold caustics in a sense that we define. In 
contrast, the walkaway and crosswell data sets both give rise to strong, nonremovable artifacts.

\end{abstract}

\maketitle

%%%%%%%%%%%%%%%%%%%%%%%%%%%%%
%%%%%%%%%%%%%%%%% section break

\section{Introduction}\label{sec intro}

In seismic acoustic imaging, borehole data refers to measurements of waves
made by receivers (sensors)  at various depths in a well;
applications include
prospecting for CO$_2$ sequestration sites or geothermal reservoirs,
and monitoring  aquifer pollution or existing hydrocarbon reservoirs.
In Vertical Seismic Profiling, the waves are excited by sources located at positions  
on the surface  \cite{BL1984,Pod};  in crosswell (or crosshole) imaging, the sources 
are in another well  \cite{D2004,D2008,A2013,S2016}.
Compared to  data resulting from traditional seismic experiments, where both the sources 
and receivers are located on the surface,
the decreased travel distance for the waves traveling to receivers in a borehole
results in less attenuation and allows the use of higher frequency waves,
potentially resulting in more sensitive and higher resolution imaging of
material parameters in the subsurface  \cite{Pod}.
\ss

The purpose of this work is to formulate a general approach  to
the analysis of borehole seismic data,
using techniques of microlocal analysis that have previously been successful for conventional
seismic data \cite{KSV98,NS1997,N00,St00,F05,F07,FG08,FG10} (and also for a variety of other imaging problems).
As in those works, here we analyze  the relation between
features in the subsurface,
in the form of singularities of the  sound speed profile, and
the resulting singularities of the data.
In some cases, the latter accurately encode the former;
in others,  imaging
artifacts arise from the data acquisition geometry,
the presence of caustics (multi-pathing) in the
subsurface, or the interaction of the two.
For several specific borehole data geometries,
we either  show that the imaging is artifact free,
or determine the location and strength of the artifacts.
\ss

In applying microlocal analysis to inverse problems, one is applying a set of tools whose theoretical foundation
is rigorously established in the high frequency limit to problems where the data is intrinsically band-limited.
Past work has shown, however, that this can be a fruitful approach,
since in practice frequencies do not have to be very high
for the high frequency limit to  be a good enough approximation
that useful  conclusions can be drawn regarding the structure and strength of
artifacts in images produced from the data.
\ss

For simplicity,  assume the Earth is  $\R^3_+$; its surface, $\R^2=\{x_3=0\}$, is flat; and  $x_3$
increases  with depth.
Throughout, we will assume that  the set $\Sigma_R$ of receivers occupies a line segment located
in a vertical borehole along the positive $x_3$-axis,
$$\Sigma_R:=\{(0,0,r): r_{min}< r< r_{max}\}.$$

 In Vertical Seismic Profiling, waves are generated by impulses located at a set $\Sigma_S$
 of sources located on the surface $\R^2$, scatter off of  features in the subsurface, and
are then measured at points of $\Sigma_R$ at all times $t\in T=(t_{min},t_{max})$.
 By contrast, for crosswell imaging  the sources  are located in another  borehole some
 distance from the  borehole containing the receivers.
Possible data sets $\D=\Sigma_S\times \Sigma_R\times T$ may be
crudely classified by their dimensionality,
depending on whether $S$ is zero-, one- or two-dimensional.
If $\dim \Sigma_S=0$, then the source set is at a single offset, or at most a discrete set of points; in this case,
$\dim \D=2$, so that the data set is underdetermined, which is not of interest for the questions we pose.
Instead, we study three basic  data sets:
\bs

\ni{\bf Overdetermined} ($\!\!\dim \D=4$):
\ms

For the {\bf dense array} data acquisition geometry (also called 3D Vertical Seismic Profiling),
$\Sigma_S\subset\partial\R^3_+\setminus{\mathbf 0}=\R^2\setminus (0,0)$  is an open subset on the surface.
\bs

\newpage

\ni{\bf Determined} ($\!\!\dim \D=3$):
\bs

In the {\bf crosswell} geometry, the sources are located in  another vertical borehole,
along a line segment,
$$\Sigma_S:=\left\{(s_0,s): s_{min}< s< s_{max}\right\},\quad s_0\ne (0,0).$$
\ms

For the {\bf walkaway} geometry, $\Sigma_S$ is contained in a line on the surface, 
passing over the borehole top
$(0,0,0)$. Without loss of generality, we assume that
$$\Sigma_S:=\left\{(s,0,0): 0<s_{min}< s< s_{max}\right\}.$$
\ms

\ni{\bf Remark.} To avoid having to deal with uninteresting degeneracies
arising from sources located very close to the top of the borehole, for all of the
geometries we  only consider sources $S=(s_0,0)=(s_1,s_1,0)$ with $ |(s_1,s_2)|>\epsilon>0$.
\ms

Our goal is to study, for each of the data sets $\D$ considered here,   the formal linearization $d\F$
about a known smooth background sound speed, $c_0(x)$
of the sound speed-to-data map, restricted to  $\D$.

In Sec. \ref{sec model}, we describe the model and  the linearization of the forward scattering
operator and  recall the Fourier integral operator theory needed in the paper.
(The  $C^\infty$ singularity theory,
describing the types of degeneracies of smooth functions
- folds, blowdowns, submersions with folds and cross caps - that
we will use to  understand the structure of the forward
and normal operators for the various data sets,
 is included as an appendix, Sec. \ref{sec sings}.)
\ss

The remainder of the paper analyzes,
for the data acquisition geometries described above,
 the linearized  scattering operator $F:=d\F$,
which is a Fourier integral operator; the geometry of its canonical relation;
and the implications for the associated normal operator, $F^*F$.
We  focus on the overdetermined dense array
data set, for which the results are the most positive.
First, in Thm. \ref{thm no caustics},
we show that if the ray geometry of the background sound speed
has no caustics, and satisfies two additional technical assumptions, then
$F$ satisfies the Traveltime Injectivity Condition,
introduced in  \cite{KSV98,NS1997}
when both the sources and receivers are on the surface.
This implies  that the normal operator
is a pseudodifferential operator, and filtered backprojection does not give rise
to artifacts in the images.
We first prove this for  a constant background sound speed  in Sec. \ref{sec dense constant},
where the additional assumptions are unnecessary.
This is followed by the analysis in Sec. \ref{sec dense nocaustic}
of  variable $c_0$ with no caustics, where  the additional assumptions are used.
\ms

Then, in Sec. \ref{sec dense fold} we study the situation for the dense array geometry
when the most commonly encountered form of caustics (or multipathing) is present.
We formulate a notion of {\it caustics  of fold type} appropriate  for this setting
and show that, in the presence of caustics no worse than this type, 
the canonical relation of $F$, while degenerate,
has  a structure, that of a {\it folded cross cap}, introduced previously 
by two of the authors in the context of marine seismic imaging \cite{FG08}; 
see Thm. \ref{thm folded cc} and Def. \ref{def folded crosscap}.
This allows a precise description of the normal operator and
characterization of  imaging artifacts microlocally away from high codimension sets.
\ms

In contrast, for the crosswell and walkaway data sets, the microlocal geometry is less favorable,
resulting in strong, nonremovable artifacts.
In Sec. \ref{crosswell} we analyze their normal operators
when the background sound speed $c_0$  is constant.
Calculating and analyzing their canonical relations  shows that
attempted inversion by filtered back projection
results in strong artifacts; see Thms. \ref{thm crosswell} and \ref{thm walkaway}.
Since these are  already badly behaved  for a constant background sound speed,
we do not pursue the analysis of the
crosswell and walkaway geometries for variable backgrounds.

%%%%%%%%%%%%%%%%%%%%%%%%%%%%%
%%%%%%%%%%%%%%%%% section break

\section{Scattering model, microlocal analysis and singularity theory}\label{sec model}

The idea of using techniques from microlocal analysis to
study linearized seismic imaging
was introduced by   Beylkin \cite{Be85},  and expanded upon by
Nolan and  Symes \cite{NS1997} and
ten Kroode, Smit and  Verdel \cite{KSV98}.
In the forward scattering problem, acoustic waves are generated
at the surface of the Earth, scatter off of features in the
subsurface, and some of the reflected waves return to the surface
to be detected by receivers, which  in these  works were also at the surface.
The goal of the  full inverse problem is to
obtain an image of the subsurface
using measurements of the pressure field at various receivers.
Due to the strong nonlinearity of  the full problem,
 works such as \cite{Be85,NS1997,KSV98} instead considered a formal linearization
of the nonlinear sound speed-to-data map, $\F$.
The linearization maps  perturbations of  a smooth
background  sound speed in the subsurface (assumed known),  to
perturbations of the resulting pressure field at receivers.
The linearized problem was  explored further by several authors for
a variety of data acquisition geometries; see, e.g.,  \cite{N00,St00,F05,F07,FG08,FG10}.
In this section, we review the scattering model and its linearization, and set down or give
references to the basic microlocal analysis  needed for the remainder of the  paper;
a summary of the requisite singularity theory  is in the appendix, Sec. \ref{sec sings}.
Since this material, in forms suitable for what is needed here, is standard, we will keep the
presentation as brief as possible.

%%%%%%%%%%%%%%%%%%%%%%%%%%%%%
%%%%%%%%%%%%%%%%% section break

\subsection{Scattering model and normal operator}\label{subsec model}

\par We now recall the scattering problem and its linearization.
Represent the Earth as $Y=\mathbb R^3_+= \{x \in \R^3, x_3 \geq 0 \}$,
consisting of isotropic material with sound speed $c(x)$,
the recovery of which is our goal.  An impulse at a source $x=s,\, t=0$,
assumed for simplicity to be a delta waveform,
creates a pressure field $p(s;x,t)$ which
solves the acoustic wave equation,
\begin{eqnarray}\label{eqn full wave}
\frac{1}{c^2(x)} \frac{\partial^2p}{\partial
t^2}(s;x,t)-\triangle p (s;x,t)&=&\delta(t)\delta(x-s)\nonumber\\
p(s;x,t)&=&0,\ \ t < 0,
\end{eqnarray}
where $\Delta$ is the Laplacian on $\R^3$.
Fixing a data acquisition set $\D=\Sigma_S\times \Sigma_R\times T$, where $\Sigma_S$ is a set of sources, 
$\Sigma_R$ is a set of receivers, and $T=(t_{min},t_{max})$ is a time interval,
the  corresponding forward map $\F=\F_{\D}$ is the sound speed-to-data map,
$c \buildrel{\F}\over{\to} p(s;r,t)|_{(s,r,t)\in\D}$;
the full inverse problem   is to reconstruct $c$ from $\F(c)$.
\ms

Due to the nonlinearity of $\mathcal F$, the works cited below  instead considered the 
formal linearization of $\F$,
assuming $c$ to be of the form $c=c_0 + \delta c$,
with $c_0 $  a known smooth  background sound speed.
Thus, associated  Green's function, i.e., background pressure field
$p_0$ satisfying \eqref{eqn full wave} for $c_0$,
is also  known (in principle). The linearization $d\F$ of $\F$  then arises from writing
 $p=p_0 + \delta p \hbox{ mod } \left(\delta c\right)^2$, where $\delta p=:(d\F)(\delta c)$ satisfies

\begin{eqnarray}\label{starstar}
\square_{c_0}(\delta p):=
\frac{1}{c_0^2(x)}\frac{\partial^2 \delta
p}{\partial  t^2}(s;x,t) -\triangle \delta p(s;x,t)&=& \frac{2}{c_0^3(x) }
\cdot \frac{\partial^2 p_0}{\partial t^2}\cdot \delta c(x)\nonumber\\
   \delta p &=&0,\ \ t < 0.
\end{eqnarray}
We denote $d\F$ by $F$, using a subscript $\D$ if needed for clarity.
Thus,
\be\label{def F}
F(\delta c)(s,r,t)=\square_{c_0}^{-1}\left(\frac{2}{c_0^3(x) }
\cdot \frac{\partial^2 p_0}{\partial t^2}\cdot \delta c(x)\right)(s;x,t)\bigg|_{x=r},
\ee
where $\square_{c_0}^{-1}$ is the forward  solution operator.
One assumes that $\delta c$ is supported at a positive distance from all sources $s$, so that its
product above with the Green's function with pole at $s$ only involves the singularities of $p_0$
on the wavefront, and not at $s$.
(The linearization can be justified in terms of Fr\'echet differentiability of $\F$ between certain
pairs of function spaces \cite{KiRi14}.)
\bs

Beylkin \cite{Be85} showed that, for a single source on the surface and an open set of
receivers, also on the surface, if caustics do not occur
for the background sound speed, then the normal operator $N:=F^*F$ is a pseudodifferential
operator $ (\Psi DO)$.
For more general ray geometries, to avoid degenerate situations one makes two assumptions:
$(i)$ no single (unbroken) ray connects a source to a
receiver; and $(ii)$ no ray originating in the subsurface grazes
$\Sigma_S$ or $\Sigma_R$. Under these assumptions, in the case of a {\it
single} source and receivers ranging over an open subset of the
surface, $\{x_3=0\}$, Rakesh \cite{Rak88} showed that $F$ is a Fourier integral
operator (FIO) in the sense of H\"ormander \cite{Hor}, and this was extended to other 
data sets $\D$ in \cite{Ha91,NS1997,KSV98}.
The assumption $(ii)$ ensures that
$F$ is an FIO, and (i)
ensures that the composition $F^*F$ makes sense on distributions.
For borehole data sets we replace these assumptions by suitably filtering (muting) the data.
\ms

The invertibility  of $F$ (modulo $C^\infty$)
was established in \cite{Be85,NS1997,KSV98}
under various combinations of assumptions on the data acquisition set $\D$
(assumed to be a smooth manifold in $\partial Y\times \partial Y$)
and the ray geometry for the background sound speed
$c_0$. In these cases, the FIO $F$ is
associated with a canonical relation $C\subset T^*\D \times T^*Y$ which satisfies the so-called
{\it traveltime injectivity condition} (TIC) described below. By the standard theory of FIO with
nondegenerate canonical relations, it follows that the normal operator
$N$ is a pseudodifferential
operator on $Y$, $N\in\Psi(Y)$; furthermore, $N$ is elliptic (and hence invertible microlocally) 
under an illumination condition.
If $Q\in \Psi(Y)$ is a left-parametrix for $N$ (i.e., a left-inverse modulo $C^\infty$), then
$Q\circ  F^*$ is a left-parametrix for $F$. This implies the injectivity of $dF$ mod $C^\infty$,
so that the singularities of $F(\delta c)$ determine the singularities of $\delta c$, as well as
giving an approximate reconstruction formula via filtered backprojection. In \cite{St00}, part
of the TIC  is relaxed, but the composition $F^*F$ is still covered by the standard clean
composition calculus for FIO.
\smallskip

On the other hand, combinations of data sets  and background ray geometries for which
the TIC is violated were studied in \cite{N00,F05,F07,FG08,FG10}. In each of these, the
composition forming the normal operator lies outside the clean intersection calculus
and $N$ is
not a  pseudodifferential operator.
The wavefront relation of $N$ is larger than that of a $\Psi$DO,
including an additional part, which gives rise to artifacts in the image when attempting filtered back projection;
the strength of that artifact depends both on the geometry of $\D$, the nature of the multi-pathing (if any), the
background ray geometry, and their interaction.
\smallskip

%%%%%%%%%%%%%%%%%%%%%%%%%%%%%
%%%%%%%%%%%%%%%%% section break

\subsection{Microlocal analysis}\label{subsec microlocal}

We now recall some basic definitions and results from  the theory of FIOs \cite{Hor}.
Let $X$ and $Y$ be smooth manifolds, of (possibly different) dimensions, $n_X,\, n_Y$, resp.
A {\it Fourier integral operator} is a continuous linear map \linebreak$A:\mathcal E'(Y)\to \mathcal D'(X)$,
whose Schwartz kernel is a locally finite sum of oscillatory integrals of the form
$$K_A(x,y)=\int_{\R^N} e^{i\varphi(x,y;\theta)} a(x,y;\theta)\, d\theta,$$
where $\varphi$ is a nondegenerate operator phase function on $X\times Y\times (\R^N\setminus 0)$, for some $N\ge 1$, and $a$ is a H\"ormander class amplitude of order $\mu$ and type $(1,0)$.
The {\it order} of $A$ is defined to be
$$m:=\mu+\frac{2N-n_X-n_Y}4,$$
and   the {\it canonical relation} of $A$ is
$$C_A:=\big\{ (x,d_x\varphi; y, -d_y\varphi)\, : \, (x,y;\theta)\in supp(a),\, d_{\theta}\varphi(x,y;\theta)=0\,\big\}\subset (T^*X\setminus 0)\times (T^*Y\setminus 0).$$

If $WF(\cdot)$ denotes the $C^\infty$ wavefront set of a distribution,
the {\it wavefront relation} of $A$,
$WF(A):=WF(K_A)'$, is the image of
the wavefront set of the Schwartz kernel of $A$ under the
map $(x,y,\xi,\eta) \to (x,\xi; y, -\eta)$; from the general theory of Fourier integral distributions,
one knows that $WF(A)\subseteq C_A$.
Thus, by the H\"ormander-Sato Lemma,  for all $u\in\mathcal E'(Y)$,
\be\label{eqn HorSat}
WF(Au)\subseteq WF(A)\circ WF(u)\subseteq C_A\circ WF(u),
\ee
where   $WF(A)$ and $ C_A$  are considered as relations from
$T^*Y\setminus 0$ to $T^*X\setminus 0$.
\medskip

For any canonical relation $C \subset
(T^*X  \setminus 0) \times (T^*Y \setminus 0)$ and $m\in\R$,
$I^m(X,Y;C)$  denotes the class of properly supported
$m$-th order FIOs $A$ with $C_A\subseteq C$. Thus, for any $A$ in this class, 
$WF(Au)\subseteq C\circ  WF(u),\, \forall u\in\mathcal E'(Y)$.
\medskip

A generalization  of Fourier integral operators are the {\it paired Lagrangian} operators
of Melrose and Uhlmann \cite{MeUh79} and Guillemin and Uhlmann \cite{GuUh81}. These are
associated to cleanly intersecting pairs of canonical relations,
$C_0,C_1\subset T^*X\times T^*Y$, and are indexed by bi-orders $(p,l)\in\R^2$.
We will not need the definitions and characterizations of these operators,
but note two properties  for later use. First, if $A\in I^{p,l}(X,Y;C_0,C_1)$,
then
\be\label{eqn iplwf}
WF(A)\subseteq C_0\cup C_1.
\ee
 Secondly, microlocally away from $C_0\cap C_1$,
\be\label{eqn ipl orders}
A\in I^{p+l}(C_0\setminus C_1)\hbox{ and }A\in I^{p}(C_1\setminus C_0).
\ee
\medskip

 Now let $C_1 \subset (T^*X
\setminus 0) \times (T^*Y \setminus 0)$ and $C_2 \subset (T^*Y
\setminus 0) \times (T^*Z \setminus 0)$ be two canonical relations,
and $A_1 \in I^{m_1}(X,Y;C_{1})$ and $A_2 \in I^{m_2}(Y,Z;C_{2})$.
If $ C_1 \times C_2 $ intersects $T^*X \times \Delta_{T^*Y}\times
T^*Z$ {\it transversely}, then H\"{o}rmander  proved that $A_1 \circ
A_2 \in I^{m_1+m_2}(X,Z;C_{1} \circ C_{2})$ where $C_{1} \circ
C_{2}$ is
the composition of $C_1$ and $C_2$ as relations in $T^*X \times T^*Y$
and $T^*Y \times T^*Z$. Duistermaat and Guillemin \cite{DuGu75} and
Weinstein \cite{We75}
extended this to the case of {\it clean
intersection} and showed that if $C_1 \times C_2$ and $T^*X \times
\Delta_{T^*Y} \times T^*Z$ intersect
cleanly with
excess $e$ then,
as in the transverse case,  $C_{1} \circ C_{2}$
is again a smooth canonical relation, and
 $A_1 \circ A_2 \in I^{m_1+m_2+e/2}(X,Z;C_1 \circ
C_2)$.
\smallskip

We say that a canonical relation $C\subset T^*X\times T^*Y$ satisfies the {\it traveltime injectivity
condition} (TIC) \cite{NS1997,KSV98}
(equivalent to the earlier {\it Bolker} condition in tomography \cite{Gu85}) if the natural projection
to the left, $\pi_L:C\to T^*X$, satisfies the following two conditions. First,
\be\label{cond tic imm}
\pi_L\hbox{ is an immersion, i.e., }d\pi_L\hbox{ is injective everywhere.}
\ee
(By results for general canonical relations, this is equivalent with $\pi_R:C\to T^*Y$ being a submersion, i.e., $d\pi_R$ is surjective.) Secondly,
\be\label{cond tic inj}
\pi_L\hbox{ is globally injective.}
\ee
(Note that \eqref{cond tic imm} already  implies that $\pi_L$ is {\it locally} injective;
\eqref{cond tic inj} demands that the injectivity holds globally.)
\medskip

If $A\in I^m(X,Y;C)$, then $A^*\in I^m(Y,X; C^t)$. If $C$ satisfies  the TIC,
then it follows from  \eqref{cond tic imm} that the composition
$A^*A$ is covered by the clean intersection calculus,
with excess $e=\dim(X)-\dim(Y)$; furthermore  \eqref{cond tic inj}
implies that $C^t\circ C\subseteq \Delta_{T^*Y}$,
the diagonal of $T^*Y\times T^*Y$. Thus,  the normal operator
\be\label{eqn clean}
N:=A^*A\in I^{2m+\frac{\dim(X)-\dim(Y)}2}(Y,Y;\Delta_{T^*Y}),
\ee
i.e., is a pseudodifferential operator on $Y$. $N$ is elliptic if $A$ is, which
in applications corresponds to an illumination condition. In that case, $N$ admits a left parametrix 
$Q\in \Psi^{-2m}(Y)$,  i.e., $Q\circ N-I$ is a smoothing operator, and then $QA^*$ is a left parametrix for $A$, 
so that, for all $u\in\mathcal E'(Y)$,  $Au \!\mod C^\infty$ determine $u \!\mod C^\infty$.
\medskip

However, in many  inverse problems,  the TIC condition  fails,
and  to understand the possibility of imaging using filtered back projection,
it is important to  analyze the composition $A^*A$ and the nature of
the resulting normal operator, $N$.
Any  component of the wavefront relation of $N$ in the complement of the diagonal
$\Delta_{T^*Y}$ will produce {\it artifacts}, i.e.,
features in $Nu$ which are not present in $u$.
It turns out that the geometry of
the canonical relation $C$,
as expressed by  degeneracies of  projections $\pi_L$ and $\pi_R$,
if they exist,
plays an important role in determining the nature, location and strength of artifacts.
\smallskip

It is known that if either
$d\pi_L$ or $d\pi_R$ has maximal rank, so does the other one and
we say that the canonical relation  $C$ is {\it nondegenerate}.
In this case the
composition  $C^t\circ C$ is clean, and \eqref{eqn clean} holds.
\smallskip

On the other hand, if  $C$ is degenerate
(the differentials of the
projections fail to be  of maximal rank),
there is no general theory that applies to the compositions $C^t\circ C$ and $A^*A$.
However, certain particular geometries  have
been analyzed, and one in particular is relevant here, for the dense array in the presence of fold caustics.
\smallskip

When one of the projections
drops rank, then the other one does, too, and their coranks are the same.
However, although corank($d \pi_L$)=corank($d \pi_R$) at all points, the
two  projections might have the same type of singularity, or quite
different ones.
The singularities needed in this article are
{\it blowdowns, folds, submersion with folds} and {\it cross caps}.
We  refer the reader to Sec. \ref{sec sings} for a concise summary of these classes;
see \cite{GoGu,Wh45,mo1} for more background,
and to Sec. \ref{sec dense fold} for the existing composition calculus \cite{FG08}, originating in marine seismic 
imaging, that we show is relevant
for the dense array with fold caustics.

%%%%%%%%%%%%%%%%%%%%%%%%%%%%%
%%%%%%%%%%%%%%%%% section break

\section{Dense array: constant $c_0$}\label{sec dense constant}

We start with the dense array geometry,
with sources  $S=(s_1, s_2, 0)$ in an open subset $\Sigma_S$  of the surface,
and receivers  in  the borehole $\Sigma_R$,  $R=(0,0,r),\, r\in (r_{min},r_{max})$.
The perturbation in sound speed is a function of the three variables, $y=(y_1,y_2,y_3)$, 
while the resulting data is a function of four variables,
$(s_1, s_2,r, t)=(s,r, t)$.
We make the following assumptions, the first of which is standard in the literature.
\ms

\basm \label{separation}
The perturbation $\delta c$ in the sound speed
has compact support at a positive depth below
the surface.
\easm
\ms

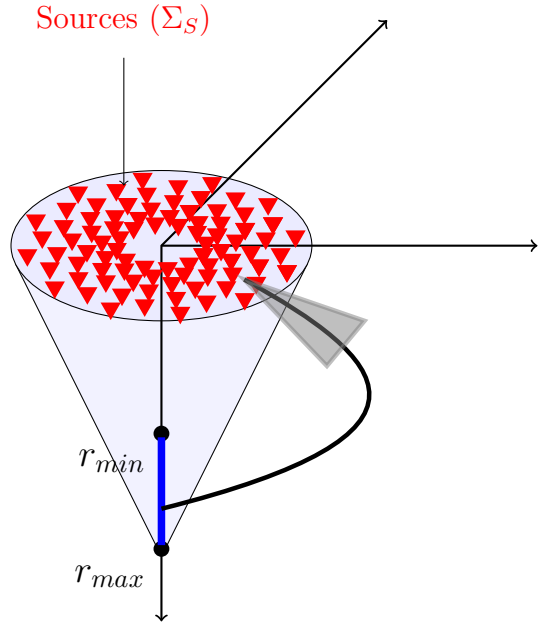
\begin{figure}
  \centering
\begin{tikzpicture}[scale=1,
  triangle/.style = {fill=red, regular polygon, regular polygon sides=3, rotate=180, scale=0.4 },
  dot/.style = {circle, fill=black, minimum size=#1,
    inner sep=0pt, outer sep=0pt},
  dot/.default = 6pt  % size of the circle diameter
  ]
  \def\a{2} \def\b{1} % Ellipse defined by p(t)=a*cos(t)+b*sin(t)
  \def\theta{194} % Value of t for tangent
  \def\vertex{1/sin(\theta)} % Vertex vertical coordinate

  % Set some parameters
  \def\freq{30}
  \def\phase{22}
  \edef\ind{1}

% Draw the cone
  \filldraw [color=blue!5!white, draw=black]  ({2*cos(\theta)},{sin(\theta)}) -- (0,{\vertex}) -- ({-2*cos(\theta)},{sin(\theta)}); % Middle point is Cone Vertex 
  \filldraw [color=blue!8!white, draw=black] (0,0) ellipse (\a cm and \b cm);

% Draw the axis
  \draw [thick, <->] (3,3) node[above right]{} -- (0,0) -- (0,-5) node[left]{}; % axes
  \draw [thick, ->] (0,0) -- (5,0) node[right]{} ; % axis

%Draw the receivers
  \foreach \radius in {.3,.4,.5,.6,.7,.8,.9}{
     \foreach \alpha in {0,\freq,...,360}{
        \draw  ({\radius*\a*cos(\alpha+(-1)^\ind*\phase)},{\radius*\b*sin(\alpha+(-1)^\ind*\phase)}) node[triangle] {};
    }
      \pgfmathparse{\ind+1}
      \xdef\ind{\pgfmathresult}
    }

% Draw Receivers
  \node[dot=6pt,label=below left:\large $r_{min}$] at (0,-2.5){};
  \node[dot=6pt,label=below left:\large $r_{max}$] at (0,{\vertex+0.1}){};
  \draw[line width=1mm,color=blue] (0,-2.55) -- (0,{\vertex+.15});

% Draw ray unbroken filtered ray
  \draw[line width=0.6mm,color=black] (0,-3.5) .. controls (4,-2.5) and (3,-1.5) .. (1.1,-0.45);
 \filldraw[line width=0.6mm,color=gray, opacity=0.5] (1.1,-.45) -- (2.7,-1) -- (2.2,-1.6) -- cycle ;

    \draw [->] (-0.5,2.5) node[yshift=0.5cm,color=red]{Sources ($\Sigma_S$)}-- (-0.5,0.8);        
    
\end{tikzpicture}
\caption{Illustration of data acquisition geometry and filtering: contributions to the data from unbroken rays such as that illustrated here are filtered out by removing data associated to those rays arriving from nearby directions, as indicated by the gray cone.}\label{nsop_fig}
\end{figure}

\basm \label{nsop}
For any unbroken ray connecting a source to a receiver which intersects the support of the
reflectivity function,  its contribution to the data has been muted by application of the filter  described
below.
(See Figure \ref{nsop_fig} for an illustration of this.)
\easm
\ms

{\bf Filter Construction:} Suppose a ray connects a source to a receiver located at
$r_0\in\Sigma_R$ and that this ray arrives at $r_0$ in a direction $\hat\rho_0$. Let $\rho_0$ be
the orthogonal projection of $\hat\rho_0$ onto $T_{r_0}\Sigma_r$.
\ms

(i) Localize the data $d(s,r,t)$
by multiplying it by a cutoff function $\chi_1(r)$ supported near $r=r_0$, and then take the partial
Fourier transform of $\chi_1d(s,r,t)$ with respect to $r$ to get $\widehat{d_1}(s,\rho,t)$, where $\rho$ is
the  Fourier variable dual to $r$.
\ms

(ii) Multiply $\widehat{d_1}$ by a cutoff function $\chi_2(\rho)$ which
is homogeneous of degree $1$ and vanishes in a conic neighborhood of the direction of
$\rho_0$ whenever $\widehat{d_1}(s,\lambda\rho_0,t)$ is {\it not} rapidly decaying as
$\lambda\to\infty$.
\ms

(iii)  Apply the inverse Fourier transform (w.r.t. $\rho$) to $\chi_2\widehat{d_1}$,
and use the result as the suitably modified data, referred to in Assumption \ref{nsop} above.
\ms

In this section and the next, we show that
the filtered linearized
scattering operator satisfies the Traveltime Injectivity Condition:

\medskip

\begin{theorem}\label{thm no caustics}
Suppose, in addition to Assumptions \ref{separation} and \ref{nsop},
 the ray geometry of a smooth background sound speed $c_0(x)$
 satisfies Assumptions \ref{simple_geom} and \ref{no_grazing} below.
 Then  the linearized   scattering operator for the dense array data set,
 $F:\mathcal E'(\R^3_+)\to \mathcal D'(\mathbb D)$
is a Fourier integral operator,
$F\in I^\frac34(\D,\R^3_+;C)$.

 If $c_0$ also has no caustics,
 then the canonical relation $C\subset T^*\D\times T^*\R^3_+$ satisfies the Traveltime Injectivity Condition
 \eqref{cond tic imm}, \eqref{cond tic inj}, and thus the normal operator is a  pseudodifferential operator of order 2,
$F^*F\in \Psi^2(\R^3_+)$.
\end{theorem}

% \documentclass[parskip]{scrartcl}
% \usepackage{tikz}
% \usetikzlibrary{shapes.geometric}

% \begin{document}

%Comment out the lines above and the \end{document} line at the end of this document once ready to incorporate the figure into main document

% \end{document}

% %%% Local Variables:
% %%% mode: latex
% %%% TeX-master: t
% %%% End:

In the current section, we consider first the model case of constant sound speed,  normalized to $c_0=1$, for
which Assumptions \ref{simple_geom} and \ref{no_grazing} hold automatically.
We will show that the canonical relation of the  linearized forward scattering operator $F$
satisfies the   traveltime  injectivity condition \eqref{cond tic imm}, \eqref{cond tic inj}.
As discussed in Sec. \ref{subsec microlocal}, this implies that the normal
operator  $F^*F$ is a pseudodifferential operator,
and a perturbation $\delta c$ of the sound speed can be reconstructed from $F(\delta c)$ 
by filtered backprojection.
\ms

\begin{proof}[Proof of Thm. 3.1]
In the case of constant background sound speed $c_0$, we   compute
the canonical relations of $F=d\F$
for each of the data geometries (dense array here; crosswell and walkaway in
Sec. \ref{crosswell}), restricting to the various data sets the basic phase function
\be\label{eqn phi basic}
\phi(s,r,t;\omega)=\left(t-|y-R|-|y-S|\right)\omega,
\ee
where $S=S(s)$, $R=R(r)$ and $y$ denote a general source, receiver and a point in the
 subsurface,  respectively, and $\omega\in\R\setminus 0$ is a phase variable. Note that the function $\phi$ 
 defined in (\ref{eqn phi basic}) is a non-degenerate phase function as can be easily verified by 
 noting that $\omega\neq 0$ and also using Assumption \ref{nsop}.
The Schwartz kernel of  $F$ is
$$K(s,r,t,y)=\int e^{i \phi(s,r,t,y; \omega)} a(s,r,t,y;\omega)\,  d\omega,$$
where $\phi$ is given by
\bast
\phi(s, r, t, y; \omega)&=&\left(t-\left|y-S\right|-\left|y-R\right|\right)\omega\\
&=&\left(t-\sqrt{(y_1-s_1)^2+(y_2-s_2)^2+y_3^2}-\sqrt{y_1^2+y_2^2+(y_3-r)^2}\right)\omega
\east
and $a\in S^2_{cl}$ is a classical symbol of order 2.
Thus, $F$ is a Fourier integral operator of order $m=2+\frac12-\frac{4+3}4=\frac34$ 
associated with the canonical
relation $C \subset T^*\R^4 \times T^*\R^3$ parametrized by $\phi$.
See \cite{NS1997} for a discussion of $F$  as an FIO  for general data sets.

Coordinates on the 7-dimensional  $C$
can be taken to be  $(s,r,y,\omega)\in\R^6\times(\R\setminus 0)$:
\bea\label{eqn C constant v1}
C&=&\Big\{\big(s_1, s_2, r, A+B, \frac{y_1-s_1}{A} \omega,  \frac{y_2-s_2}{A} \omega, 
\frac{y_3-r}{B} \omega, \omega; \nonumber \\
 & &  \quad  y_1, y_2, y_3, \left(\frac{y_1-s_1}{A}+  \frac{y_1}{B}\right)\omega,  
 \left(\frac{y_2-s_2}{A}+  \frac{y_2}{B}\right)\omega,
 \left(\frac{y_3}{A}+  \frac{y_3-r}{B}\right)\omega\big) \Big\},
 \eea

 \ni where $A=\sqrt{(y_1-s_1)^2+(y_2-s_2)^2+y_3^2}$  \, and $B=\sqrt{y_1^2+y_2^2+(y_3-r)^2}$.
Thus, $C \subset \left(T^*\R^4\setminus 0\right) \times \left(T^*\R^3\setminus 0\right)$,
and we see that the projections to left and right are
\be\label{eqn piL}
\pi_L(s, r, y, \omega)=\left(s_1, s_2, r, A+B, \frac{y_1-s_1}{A} \omega,  \frac{y_2-s_2}{A} \omega, 
\frac{y_3-r}{B} \omega, \omega                                                                         \right)
\ee
 and
 \be\label{eqn piR}
 \pi_R( s, r, y,\omega)=\left( y_1, y_2, y_3, \left(\frac{y_1-s_1}{A}+  \frac{y_1}{B}\right)\omega,  
 \left(\frac{y_2-s_2}{A}+  \frac{y_2}{B}\right)\omega, \left(\frac{y_3}{A}
 +  \frac{y_3-r}{B}\right)\omega                              \right).
 \ee
Since $\pi_R$ has identity in $y$ variables, we have that  rank $d\pi_R=3+\hbox{ rank}\left(\frac{D\eta}{D(s, r, \omega)}\right)$,
where $\eta$ is dual to $y$,
and the minor $\frac{D\eta}{D( s, \omega)} $ is
 $$\frac{D \eta}{D(s, \omega)}=\left[ \begin{array} {ccc}  -\frac{(y_2-s_2)^2+y_3^2}{A^3} \omega 
 & \frac{(y_1-s_1)(y_2-s_2)}{A^3}\omega & \frac{y_1-s_1}{A}+  \frac{y_1}{B} \\
\frac{(y_1-s_1)(y_2-s_2)}{A^3}\omega & -\frac{(y_1-s_1)^2+y_3^2}{A^3}\omega 
& \frac{y_2-s_2}{A}+ \frac{y_2}{B}\\
\frac{y_3(y_1-s_1)}{A^3}\omega & \frac{(y_2-s_2)y_3}{A^3}\omega & \frac{y_3}{A}+ \frac{y_3-r}{B}\\
\end{array} \right],$$
\smallskip

\ni  which has determinant
\begin{eqnarray*}
 { \omega^2y_3}{A^{-3}}\left(1+ \frac{(y-S) \cdot (y-R)}{AB} \right).
\end{eqnarray*}
Using $\omega\ne 0,\,  y_3 > 0$,   the Cauchy-Schwarz   inequality and Assumption  \ref{nsop},
one sees that $\hbox{det}\left[\frac{D \eta}{D( s, \omega)}\right] \ne 0$,
so that $\hbox{rank}(d\pi_R)=6$.
It follows that  $d\pi_R $ has maximal rank and $\pi_R$ is a submersion;
hence    $\pi_L$  is an immersion, and $C$ is a nondegenerate canonical relation.
\ms

To verify the Traveltime  Injectivity Condition, it remains to show the injectivity of  $\pi_L$;
we do this using Assumption \ref{nsop}.
The unit vectors $\left(y-\left(0,0,r \right)\right)/B$ and
$\left(y-\left(s_1,s_2,0\right)\right)/A$  point
  to $y$  from the source $S=(s_1, s_2, 0)$ and  from the receiver $R=(0,0,r)$, resp.
  In terms of these, the condition in Assumption \ref{nsop} is that
 $$\frac1{A}\left(y-\left(s_1,s_2,0\right)\right) +   \frac1{B}\left(y-\left(0,0,r \right)\right)  \ne0.$$

To prove that $\pi_L$ is injective, let us consider
$$S=(s_1,s_2,0);\quad R=(0,0,r);\quad y = (y_1,y_2,y_3);\quad \omega;$$
and
$$\tilde{S} = (\ts_1,\ts_2,0);\quad \tilde{R} = (0,0,\tr); \quad \ty = (\ty_1, \ty_2, \ty_3);
\quad \tilde\omega;$$
such that $\pi_L(s,r,y,\omega)=\pi_L(\ts,\tr,\ty,\ww)$. Then

\begin{eqnarray}\label{piL injectivity equalities}
s_1 &=& \ts_1,\quad s_2 = \ts_2,\quad r = \tilde{r},\quad  \omega  = \ww,
\quad A+B = \tilde{A} +\tilde{B}\\
& & \nonumber \\
\frac{y_1-s_1}{A} \omega &=&  \frac{\ty_1-\ts_1}{\tilde{A}} \ww,
\quad  \frac{y_2-s_2}{A} \omega =  \frac{\ty_2-\ts_2}{\tilde{A}} \ww, \quad
 \frac{y_3-r}{B} \omega =  \frac{\ty_3-\tr}{\tilde{B}} \ww, \nonumber\\
 \nonumber
\end{eqnarray}
where

\begin{equation}
A=|y-S|;\quad B = |y-R|;\quad
\tilde{A}=|\ty-\tilde{S}|;\quad \tilde{B} = |\ty-\tilde{R}|.
\end{equation}

The first four equalities  in \eqref{piL injectivity equalities} imply that
$\tilde S=S,\, \tilde R=R,\, \ww=\omega$,
so to prove injectivity of $\pi_L$ we only need to verify that $\ty=y$.
Defining

\begin{equation}\label{rho and sigma}
 \sigma:=\frac{y-S}{A};\qquad \tsig= \frac{\ty - \tilde{S}}{\tilde{A}};\qquad \rho := \frac{y-R}{B},
 \qquad \tilde\rho := \frac{\ty - \tilde{R}}{\tilde{B}},
\end{equation}

the last three equalities in \eqref{piL injectivity equalities}
(together with the fact that $\omega = \ww$) imply that

\begin{eqnarray}
\sigma_i=\tilde \sigma_i,\ i=1,2;\\
\rho_3=\tilde\rho_3. \label{eq_rho}
\end{eqnarray}

Recalling that $\sigma=(\sigma_1,\sigma_2,\sigma_3)$ and $\tsig = (\tsig_1,\tsig_2,\tsig_3)$ 
are unit vectors, we obtain
\begin{eqnarray*}
  \sigma_3^2=1-\sigma_1^2-\sigma_2^2=1-{\tsig}_1^2-{\tsig}_2^2={\tsig}_3^2  \\
  \Rightarrow \tsig_3=\pm\sigma_3.
\end{eqnarray*}

However,  $y_3,\ty_3>0 \Rightarrow \sigma_3,\tsig_3>0$, hence $\sigma=\tsig$, that is

\begin{eqnarray}\label{ratio}
  \frac{\ty-S}{|\ty-S|}=\frac{y-S}{|y-S|} \label{unit_eq} \ ,
\end{eqnarray}
which shows that $y$ and $\ty$ lie on the same ray emanating from $S$. 
The last equality in \eqref{piL injectivity equalities} can be expressed as
\begin{eqnarray}
   T(s,r,y)=T(s,r,\tilde y), \quad \mbox{where}\ T(s,r,y)=A(y,s)+B(y,r).
\end{eqnarray}
In the argument below, if $z\in\mathbb R^3$ and $z\neq 0$, we denote the corresponding unit vector as 
$\hat z := z/|z|$. We also denote a point on the line segment between $y$ and $\tilde y$ as 
$y_t:=y+t(\tilde y-y)$, where $t\in [0,1]$.

To prove that $\pi_L$ is an injection, we argue by contradiction. Assume that $\tilde y\neq y$;
the Mean Value Theorem then implies $\exists c\in (0,1)$ such that
\begin{eqnarray*}
  \left. \frac{d}{dt}\left\{ T(s,r,y_t) \right\} \right|_{t=c}=0\\
  \Rightarrow \nabla_yT(s,r,y_c)\cdot(\tilde y-y)=0 \\
  \Rightarrow \left(\widehat{y_c-S}+\widehat{y_c-R}\right)\cdot (\tilde y-y)=0\\
  \Rightarrow \left(\widehat{y_c-S}+\widehat{y_c-R}\right)\cdot \widehat{(\tilde y-y)}=0\\
  \Rightarrow \left(\widehat{y_c-S}+\widehat{y_c-R}\right)\cdot \widehat{(y_c-S)}=0\\
  \Rightarrow \widehat{(y_c-R)}\cdot \widehat{(y_c-S)}=-1,
\end{eqnarray*}
where we have used the fact that $\widehat{(\tilde y-y)}=\widehat{(y_c-S)}$, which in turn is true because 
$s,y,\tilde y,y_c$ all lie on a single ray emanating from $S$.
However, the last equality contradicts Assumption (\ref{nsop}), and therefore $\tilde y=y$.
\ss

Hence, $C$ satisfies the Traveltime Injectivity Condition.
Taking the clean intersection calculus (\ref{eqn clean}) into account,
since $dim(\D)-dim(\R^3_+)=1$ we should write the order of $F$ as 
$\frac34=1-\frac14$; since the `effective' order of $F$
 is 1, the composition $F^*F$ is  a pseudodifferential operator of order 2 on $\R^3_+$,
which will be elliptic at those points in $T^*\R^3_+$ where $F$ is.
Thus, under an illumination assumption,
a perturbation $\delta c$ of the sound speed can be reconstructed without artifacts
by filtered backprojection: $F(\delta c)\mod C^\infty$ determines $\delta c \mod C^\infty$.

%%%%%%%%%%%%%%%%%%%%%%%%%%%%%
%%%%%%%%%%%%%%%%% section break

 \section{Dense Array: No Caustics}
 \label{sec dense nocaustic}

In this section, we continue the proof of Thm. \ref{thm no caustics}, modifying the argument of
the previous section to a variable background speed, $c_0=c_0(x)$, satisfying
 assumptions which we  now describe.

Parametrize the maximally defined characteristic curve (i.e., a ray) departing $z$ in the direction
 \begin{eqnarray} \label{nu}
   \nu(\theta,\varphi) =  (\sin(\varphi)\cos(\theta),\sin(\varphi)\sin(\theta),\cos(\varphi))
 \end{eqnarray}
 by a smooth function
 \begin{eqnarray} \label{ray}
    \R\supseteq I \owns p \mapsto x(p;\theta,\varphi,z)\in \R^3_+,
 \end{eqnarray}
 where $x(0;\theta,\varphi,z)=z$.
 The angle $\varphi$ is the polar angle with respect to the $x_3$-axis. If the take-off angle corresponds to
 $\varphi=0,\pi$, then we change coordinates so $\varphi$ is the polar angle with respect to another axis and 
 adjust (\ref{nu}) correspondingly. In the following discussion, we proceed as though polar angles $\varphi$ are 
 the polar angles with respect to the $x_3$-axis but none of the arguments depend on this and will work just as 
 well if $\varphi$ is another polar angle. Following \cite{WWS2009}, we make the following assumption:
\basm [No Caustics] \label{simple_geom}
Assume  $sing\ supp\ (V)$ is contained in a region $\Omega\subset \R^3_+$ which is completely illuminated by 
each source and receiver with a unique minimal traveltime ray connecting each point $y\in\Omega$ 
with each $z\in\Sigma_S\cup\Sigma_R$. Also assume that there are no caustic points in $\Omega$ 
on rays issuing from any $z\in\Sigma_S\cup\Sigma_R$.
\easm

Under this assumption, we have a well-defined and smooth traveltime function, $t_{c_0}(z,y)=t_{c_0}(y,z)$,
which is the minimal travel time between any $z\in\Sigma_S\cup\Sigma_R$ and $y\in\Omega$.
We will also make the standard no-grazing ray assumption:
\basm [No grazing rays] \label{no_grazing}
Whenever $z \in sing\ supp\ (V)\subset \R^3_+$, we assume
\begin{eqnarray}
  \pder{x_3}{p}(p;\theta,\varphi,z)\neq 0,\ \mbox{when}\ x(p,\theta,\varphi,z)\in\Sigma_S,
\end{eqnarray}
\easm
\noindent which means that there are no rays emanating from the subsurface to graze $\Sigma_S$.
\begin{remark}
With this setup, the phase function $\phi$ from the previous section is replaced, as in   \cite{NS1997}, by
 \begin{eqnarray}\label{inhom_phase}
   \phi(s,r,t,y,\omega):=\omega \left(t- T\left(s,y,r\right)\right),\nonumber \\ \mbox{where}\
   T(s,y,r) := t_{c_0}\left(\left(s_1,s_2,0\right),y\right)
   +t_{c_0}\left(\left(0,0,r\right),y\right), \nonumber
\end{eqnarray}
the sum of the travel times of the incident and reflected rays, is the total travel time. 
One easily verifies, using Assumption \ref{nsop},
that  $\phi$ is a non-degenerate phase function, so that $F$ is an FIO.
\end{remark}
\smallskip

Writing $z=(z_1,z_2,z_3)=(z',z_3)$, the wavefront relation of $F$ is
contained in the canonical relation
 \begin{eqnarray} \label{reln_no_caustic}
   C= \Big\{ \big( s, r,T(s,r,y) , \omega\sigma(s,y),\omega\rho(r,y), \omega\ ;\
    y, -\omega \left[ \eta_s(s,y)+\eta_r(r,y) \right] \big)\ \big| \nonumber \\
    (0,0,r)\in \Sigma_R,\, (s_1,s_2,0)\in\Sigma_S,\,\ \omega\in\R\setminus 0  \Big\}, 
 \end{eqnarray}
where
 \bea
 \sigma(s,y):=\nabla_{z'}\tau((s_1,s_2,0),y);\ \rho(r,y):=\pder{\tau}{z_3}((0,0,r),y);\\
 \eta_s(s,y):=\nabla_y\tau((s_1,s_2,0),y);\ \eta_r(r,y):=\nabla_y\tau((0,0,r),y).
 \eea

Since $F$ is an FIO, $C'$ is a $7$-dimensional conic Lagrangian submanifold of $T^*\R^7$. We note that the
above canonical relation avoids the zero section due to $\phi$ being a Hormander-type non-degenerate phase
function. One can  check that Assumptions \ref{simple_geom}, \ref{no_grazing}  imply that we  may
parametrize $C$ using coordinates
$(r,t_{ref},\check\varphi,\check\theta,\hat\varphi,\hat\theta,\omega)$, \lb defined as follows, describing each
broken ray backwards:
a ray, departing $(0,0,r)$ in  direction
$\left(\sin(\hat\varphi)\cos(\hat\theta),\sin(\hat\varphi)\sin(\hat\theta),\cos(\hat\varphi)\right)$
and traveling for  time $t_{ref}>0$, arrives at location $y:= x(t_{ref};\hat\varphi,\hat\theta,r) \in \mathbb R^3_+
$ (see (\ref{ray}) for a reminder of the definition of the function $x$ here).
Then a ray leaving $y$ in the direction
$\left(\sin(\check\varphi)\cos(\check\theta),\sin(\check\varphi)\sin(\check\theta),\cos(\check\varphi)\right)$
arrives at $(s,0) \in \Sigma_S$, where $(s,0)=x(t_{inc};\check\varphi,\check\vartheta,y)$; here, the
travel time function $t_{inc}:=t_{inc}(y,\varphi,\vartheta)$ is the travel time needed for this ray to reach
$\Sigma_s$. Note that  $t_{inc}(y,\varphi,\vartheta)$ is smooth and guaranteed to exist by the
implicit function and the non-grazing ray assumption. Finally $t:=t_{inc}+t_{ref}$ is the two-way traveltime from
$(s,0)$ to $y$ and from $y$ to $(0,0,r)$.

% \documentclass[parskip]{scrartcl}
% \usepackage{tikz}
% \usetikzlibrary{shapes.geometric}

% \begin{document}

%%% Comment out the lines above and the \end{document} line at the end of this document once ready to incorporate the figure into main document

\begin{figure} \label{coords_fig}
  \centering
\begin{tikzpicture}[scale=1,
  triangle/.style = {fill=red, regular polygon, regular polygon sides=3, rotate=180, scale=0.4 },
  dot/.style = {circle, fill=black, minimum size=#1,
    inner sep=0pt, outer sep=0pt},
  dot/.default = 6pt  % size of the circle diameter
  ]
  \def\a{2} \def\b{1} % Ellipse defined by p(t)=a*cos(t)+b*sin(t)
  \def\theta{194} % Value of t for tangent
  \def\vertex{1/sin(\theta)} % Vertex vertical coordinate

  % Set some parameters
  \def\freq{30}
  \def\phase{22}
  \edef\ind{1}

% Draw the cone
  \filldraw [color=blue!5!white, draw=black]  ({2*cos(\theta)},{sin(\theta)}) -- (0,{\vertex}) -- ({-2*cos(\theta)},{sin(\theta)}); % Middle point is Cone Vertex 
  \filldraw [color=blue!8!white, draw=black] (0,0) ellipse (\a cm and \b cm);

% Draw the axis
  \draw [thick, <->] (3,3) node[above right]{} -- (0,0) -- (0,-5) node[left]{}; % axes
  \draw [thick, ->] (0,0) -- (5,0) node[right]{} ; % axis

%Draw the receivers
  \foreach \radius in {.2,.3,.4,.5,.6,.7,.8,.9}{
     \foreach \alpha in {0,\freq,...,360}{
        \draw  ({\radius*\a*cos(\alpha+(-1)^\ind*\phase)},{\radius*\b*sin(\alpha+(-1)^\ind*\phase)}) node[triangle] {};
    }
      \pgfmathparse{\ind+1}
      \xdef\ind{\pgfmathresult}
    }

% Draw Receivers
  \node[dot=6pt,label=below left:\large $\textcolor{blue}{r_{min}}$] at (0,-2.5){};
  \node[dot=6pt,label=below left:\large $\textcolor{blue}{r_{max}}$] at (0,{\vertex+0.1}){};
  \draw[line width=1mm,color=blue] (0,-2.55) -- (0,{\vertex+.15});

% Draw ray and transverse plane
  \node[dot=6pt,label=below left:\large $y_0$] at (4,-3.5){};
  \node[dot=6pt,label=right:\large $\textcolor{blue}{(s_0 ,0)}$] at (1.7,-.5){};  
  \draw[line width=0.2mm,color=black] (3.1,-4.8) -- (5.5,-3.5) -- (4.2,-2.2) -- (2.1,-3.3) -- cycle;
  \draw[line width=0.6mm,color=black] (4,-3.5) .. controls (3,-1.5) .. (1.7,-0.45);
  \draw[line width=0.5mm,color=orange,->] (4,-3.5) -- (6.4,-2.2) node[right]{$y_3$};
  \draw[line width=0.5mm,color=orange,->] (4,-3.5) -- (5.1,-5.2) node[right]{$y_2$};
  \draw[line width=0.5mm,color=orange,->] (4,-3.5) -- (3,-1.2) node[above right]{$y_1$};  
    %Label Parts

\end{tikzpicture}
\caption{Construction of the $y$-coordinates with the $y_1$-direction being tangent to the ray connecting $y$ to a source $(s,0)\in \Sigma_S$.}\label{coords_fig}
\end{figure}
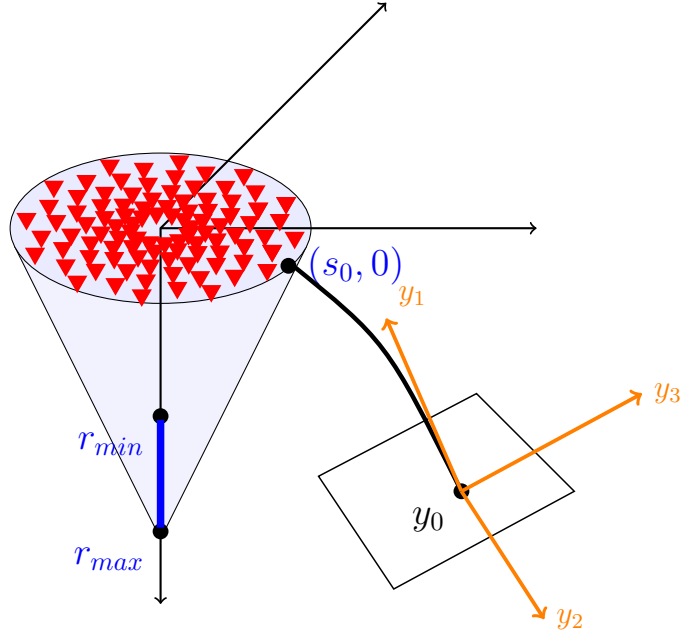

% \end{document}

% %%% Local Variables:
% %%% mode: latex
% %%% TeX-master: t
% %%% End:

We now verify that $C$ satisfies the Traveltime Injectivity Condition, starting by showing that $\pi_L$ is an 
immersion.  We check this by showing that
\begin{eqnarray} \label{j1}
  \left| \pder{(s,\sigma,T)}{(t_{ref},\hat\varphi,\hat\theta,\check\varphi,\check\theta)}  \right| \neq 0 \ .
\end{eqnarray}
Observe that the assumption of no caustics (\ref{simple_geom}) implies
\begin{eqnarray}
  \left|
  \pder{y}{(t_{ref},\hat\theta,\hat\varphi)}
  \right| \neq 0; \label{nsc}
\end{eqnarray}
 using the chain rule, it will follow that $\pi_L$ is an immersion once we establish
\begin{eqnarray}
  \left|
  \pder{(s,\sigma,T)}{(y,\check\theta,\check\varphi)}
  \right| \neq 0 \label{canon}
\end{eqnarray}
as follows. Make a change of $y$-coordinates so that the
$y_1$-direction is parallel with the velocity of a specific ray departing $y_0\in \mathbb R^3_+$ and arriving at 
$(s_0,0)\in\Sigma_S$. Let $(y_2,y_3)$ be coordinates on the plane that contains $y_0$ 
and is also orthogonal to the $y_1$-direction, as illustrated in Figure \ref{coords_fig}. 
With this choice of coordinates, the no-grazing ray and no caustics assumptions imply that
\begin{eqnarray}
\left| \pder{(s_1,s_2)}{(\check\varphi,\check\theta)}  \right| \neq 0 \label{var1}
\end{eqnarray}
and also, the no-grazing ray assumption guarantees that
\begin{eqnarray}
\left| \pder{(\sigma_1,\sigma_2)}{(y_2,y_3)} \right| \neq 0 \ . \label{var2}
\end{eqnarray}
Let $\gamma:(-\epsilon,\epsilon)\to \mathbb R^3_+, \gamma(t)=y(t)$ be a parametrisation of an open interval of the ray connecting $y_0$ to $s_0$, for a suitably small $\epsilon\in\mathbb R_+$, with $\gamma(0)=y_0$.
By construction, the $y_1$-direction is tangent to the ray connecting $y_0$ to $(s_0,0)$ and since $\left(s(\gamma(t)),\sigma(\gamma(t))\right)=(s_0,\sigma_0), \forall t\in (-\epsilon,\epsilon)$, we therefore have
\begin{equation}
\frac{1}{\dot\gamma(0)}\ \frac{d}{dt}{\left(s(\gamma(t)),\sigma(\gamma(t))\right)} = \left( \pder{s}{y_1}(y_0),\pder{\sigma}{y_1}(y_0) \right)= (0,0) \ . \label{var3}
\end{equation}
We also have
\begin{align}
\pder{T}{y_1}(s,r,y) & = \pder{t_{c_0}((0,0,r),y)}{y_1} + \pder{t_{c_0}((s_1,s_2,0),y)}{y_1} \nonumber \\
& =\pder{t_{c_0}((0,0,r),y)}{y_1}-c^{-1}(y) \ , \label{var4}
\end{align}
with the latter equality following from the definition of the $y_1$-direction. Additionally,
\begin{eqnarray}
\left|\pder{t_{c_0}((0,0,r),y)}{y_1}\right|\leq \left|\nabla t_{c_0}((0,0,r),y))\right|=c_0^{-1}(y),
\label{var5} \end{eqnarray}
where we have used the eikonal equation in the last equality. Furthermore, equality is attained in the left side of (\ref{var5}) if and only if we have scattering over $\pi$, which is ruled by Assumption \ref{nsop}. It now follows from (\ref{var4}-\ref{var5}) that
\begin{equation}
\pder{T}{y_1}(s_0,r_0,y_0) \neq 0 \ , \label{var6}
\end{equation}
for any $r_0\in (r_{min},r_{max})$.  Therefore, (\ref{var1}-\ref{var3},\ref{var6}) establish (\ref{canon}) and so we have shown that $\pi_L$ is an immersion.
\ms

We now verify that $\pi_L$ is injective. To prove this, it will be convenient to use $(s,r,y,\omega)$ as coordinates
on $C$. Suppose that $\pi_L(s,r,y,\omega)=\pi_L(\tilde s,\tilde r,\tilde y,\tilde\omega)$. Then we immediately have
$s=\tilde s,r=\tilde r, \omega=\tilde\omega$ and we deduce that
\be \label{inj_3}
\nabla_s t_{c_0}((s,0),y)=\nabla_s t_{c_0}((s,0),\tilde y) =:(\sigma_1,\sigma_2), \quad
T(s,r,y)=T(s,r,\tilde y).
\ee
To prove injectivity of $\pi_L$, it remains to show that $\tilde{y}=y$.
Condition (\ref{inj_3}) implies $y$ and $\tilde y$ lie on a common ray issuing from $(s,0)$ in the direction
\begin{eqnarray}\nonumber
\left(\sin(\varphi_s)\sin(\theta_s),\sin(\varphi_s)\cos(\theta_s),\cos(\varphi_s)\right) := \nonumber \\
\left(-\sigma_1,-\sigma_2,-\sqrt{c_0^{-2}(r,0)-\sigma_1^2-\sigma_2^2}\right)
\end{eqnarray}
for some angles $(\varphi_s,\theta_s)$. Let $p_1,p_2$ be the values of $p$ that satisfy
\begin{eqnarray*}
  y=x(p_1;\theta_s,\varphi_s,(s,0));\  \tilde y= x(p_2;\theta_s,\varphi_s,(s,0)),
\end{eqnarray*}
which then implies that
\begin{eqnarray}
T(s,r, x(p_1;\theta_s,\varphi_s,(s,0)) ) = T(s,r, x(p_2;\theta_s,\varphi_s,(s,0)) )\ . \label{varTT}
\end{eqnarray}
\ms

We can now use the same argument used at the end of Sec. 3 to show that, if we assume that $\tilde y\neq y$, 
then (\ref{varTT}) contradicts Assumption \ref{nsop}. So, under the above assumptions, $\pi_L$ is injective;
combining this with $\pi_L$ being an immersion, established earlier,
we have shown that the Traveltime Injectivity Condition is satisfied.
Thus, as in the case
of constant $c_0$,  $F^*F\in \Psi^2(\R^3_+)$, concluding the proof of Thm. \ref{thm no caustics}. \end{proof}

%%%%%%%%%%%%%%%%%%%%%%%%%%%%%
%%%%%%%%%%%%%%%%% section break

\newpage

\section{Dense Array: Fold Caustics}\label{sec dense fold}

\par
We start by formulating a notion of what it means for the ray geometry of $c_0$ to have fold caustics with
respect to  borehole data acquisition.

\subsection{Fold caustics: single receiver}\label{subsec single}
First recall the concept  for conventional seismic data, where both sources and receivers are on the surface
\cite{NS1997,N00}.
For a single receiver, $r$,
the cotangent space $\Lambda_r:=T^*_r\R^3$ is a Lagrangian submanifold
of $T^*\R^3$, on which the canonical dual variables $q=(q_1,q_2,q_3)$ are coordinates. 
The  exponential map, $\chi:=\exp\left(H_{c_0}\right):T^*\R^3\to T^*\R^3$,
of the  Hamiltonian  vector  field of the  (nonhomogeneous) symbol $\frac{1}{2}(c_0(x)^{-2}-|\xi|^2)$
 is a    (nonhomogeneous) canonical transformation of $T^*\R^3$.
Thus, the image $\Lambda_r^{c_0}:=\chi\left(\Lambda_r\right)$  is also
a Lagrangian, which is  not conic since  $\chi$ is not homogeneous in $\xi$; 
on $\Lambda_r^{c_0}$, the pushforwards $\tilde{q}:=\chi_*(q)$ by $\chi$ are coordinates.
 On $\Lambda_r^{c_0}$ there is a well-defined acoustical distance function, which is the integral of 
$c_0^{-1}$ along each bicharacteristic.

A {\it caustic} of $c_0$ (with respect to $r$) is a point $\lambda_0=(x^0,\xi^0)\in\Lambda_r^{c_0}$
where the spatial projection
$\pi_X:\Lambda_r^{c_0}\to \R^3$ has a noninvertible differential, 
and  $\lambda_0$ is a {\it fold} caustic if $\pi_X$
has a Whitney fold singularity at $\lambda_0$ (see Def. \ref{def fold}).
At such a point, $d\pi_X\left(T_{\lambda_0}\Lambda_r^{c_0}\right)$ is a hyperplane  $\Pi\subset T_{x^0}\R^3$.
For the following, assume that $\Pi$ is not vertical; otherwise, the discussion needs to be slightly
modified.
Since $dim(\Pi)=2$ and is nonvertical, $x_1, x_2$ have linearly independent gradients and thus are independent
functions on $\Lambda_r^{c_0}$. 
By Darboux's Theorem, these 
may be augmented with
$p_3:=\xi_3|_{\Lambda_r^{c_0}}$ to obtain a coordinate system on $\Lambda_r^{c_0}$ near $\lambda_0$.
On  $\Lambda_r^{c_0}$, the restrictions of the other canonical coordinates on $T^*\R^3$
are functions of $(x_1,x_2,p_3)$: $x_3= f(x_1, x_2, p_3)$ and
$(p_1, p_2) :=(\xi_1,\xi_2)|_{\Lambda_r^{c_0}}= (g_1(x_1, x_2, p_3), g_2(x_1, x_2, p_3))$.
The fold caustic  at $\lambda_0$  then implies that 
{\be\label{eqn causticfold conds}
\frac{\partial f}{\partial p_3}=0,\quad
\frac{\partial^2 f}{\partial p_3^2} \neq 0.
\ee}
Note also that the acoustical distance function
described above  is a smooth function of  $(x_1,x_2,p_3)$, since they are coordinates on $\Lambda_r^{c_0}$.

%%%%%%%%%%%%%

\subsection{Fold caustics: borehole data}\label{subsec borehole}

Now let $\D$ be the dense array data set as in the previous two sections,  for which
$\Sigma_S\subset\partial\R^3_+\setminus {\mathbf 0}\simeq \R^2\setminus (0,0)$  is an open subset and
$\Sigma_R=\left\{(0, 0, r):\, r_{min}<r<r_{max}\right\}$.
For each value of $r$, one can repeat the above constructions and analysis.
Since each $\Lambda_r=T^*_r\R^3$ is  Lagrangian, it follows that
$\Gamma:=\bigcup_{r} \Lambda_r$ is a 4-dimensional coisotropic  (or involutive) submanifold of $T^*\R^3$,
and $\Gamma$ is foliated by the family of $\Lambda_r$.
Thus, with the canonical transformation $\chi$ as in Sec. \ref{subsec single}, 
the image $\Gamma^{c_0}:=\chi\left(\Gamma\right)$ is also
a (nonconic)  four-dimensional  coisotropic submanifold of $T^*\R^3$,  
foliated by  $\left\{\Lambda_r^{c_0}:\, r_{min}<r<r_{max}\right\}$. 
At  a regular point of the spatial projection,  $\pi_X:\Gamma^{c_0}\to\R^3$,
 $\rank(d\pi_X)=3$ is maximal,
 while a caustic is a $\gamma_0$ where $\rank\left(d\pi_X(\gamma_0\right)\le 2$.
 We will demand that $\pi_X$ has at most fold singularities. Due to the difference in
 dimensions, this means  that it is a submersion with folds
 (see Def. \ref{def SWfold}).
  \medskip

\begin{defn}\label{def fold caustic}
We say that the ray geometry of $c_0$ has
\emph{at most fold caustics with respect to the borehole} $\Sigma_R$  if

(i) for each $r_{min}<r<r_{max}$,  the  only singularities of the spatial projection \lb
$\pi_X:\Lambda_r^{c_0}\to\R^3_+$ are Whitney folds; and

(ii)  the  only singularities of
$\pi_X: \Gamma^{c_0} \rightarrow \R^3_+$ are submersions with folds.
\end{defn}

\ni{\bf Remark:} One can compare  conditions (i) and (ii).
These are in fact independent of each other: The receiver-by-receiver Whitney fold condition (i) does not imply 
(ii), since the latter requires the invariantly defined Hessian to have rank two, which cannot be derived from the 
rank one Hessian coming from a Whitney fold. Conversely, (ii)  
only implies that the $\pi_X:\Lambda_r^{c_0}\to\R^3_+$ are Whitney folds
under a  tangent space condition which, while generic,   does not appear to be physically required.
This is because the restriction of a  submersion with folds  to a submanifold passing through the critical set
is not necessarily a fold; for example, the function $F(x_1,x_2)=x_1x_2$ is a submersion with folds 
$\R^2\to \R^1$, but restricted to either axis it is not a Whitney fold $\R\to\R$.

We also mention that this Def.  \ref{def fold caustic} is related to but differs from the notion of fold caustic
formulated in  \cite{FG08}  for another overdetermined data set, the  marine data acquisition geometry.

The main result of this section is the following; the  terminology used and the consequences 
for imaging are explained after its statement.
\ms

\begin{thm}\label{thm folded cc}
Under the fold caustic assumption,  and a small  slope assumption on the caustic surface 
(see \eqref{cond dip}),
the linearized forward operator $F$ is a Fourier integral operator, $F\in I^{\frac34}(\D,\R^3_+;C)$,
 whose canonical relation $C$
is a folded crosscap in the sense of Def. \ref{def folded crosscap} below, 
away from a possible set of codimension at least four.
\end{thm}

We now recall the class of
degenerate canonical relations and Fourier integral operators referred to in the theorem, which was
originally introduced by two of the authors in the context of marine seismic imaging.
Suppose that $\dim(X)=n+1$, $\dim(Y)=n$ and $C\subset \left(T^*X\setminus 0\right) 
\times \left(T^*Y\setminus 0\right)$ is a canonical relation, so that
$$\dim(T^*Y)=2n<\dim(C)=2n+1<\dim(T^*X)=2n+2.$$

\begin{defn}\label{def folded crosscap}
  {\bf \cite{FG08}}  The canonical relation
$C$ is a \emph{folded cross cap} if
\smallskip

(i) $\pi_R:C\to T^*Y$ is a submersion with folds (see
Def. \ref{def SWfold}) and the image of its critical manifold, $\pi_R\left(\Sigma\left(\pi_R\right)\right)$, is a
nonradial hypersurface in $T^*Y$;
\smallskip

(ii)  $\pi_L:C\to T^*X$  is a cross cap (see Def. \ref{def crosscap}) and
%(iii)
$\pi_L\left(\Sigma\left(\pi_L\right)\right)$, which is   a codimension three,
 immersed submanifold in $T^*X$, is also nonradial.
\end{defn}

\begin{remark}\label{rem fcc def}
Recall that \emph{nonradial} means that the restriction of the canonical 1-form does not vanish anywhere. 
Also, from \cite{FG08}
one knows that $\pi_L\left(\Sigma\left(\pi_L\right)\right)$ must be maximally noninvolutive, i.e., the 
restriction to it of the canonical two form   has maximal possible rank everywhere, 
which on the $(2n-1)-$ dimensional 
$\pi_L\left(\Sigma\left(\pi_L\right)\right)$ is $2n-2$.
\end{remark}

For  a folded cross cap, the composition $C^t\circ C$ lies outside the clean intersection calculus,
but the following holds.

\begin{thm}\label{thm folded crosscap}   {\bf \cite{FG08}} If $C$ is a folded cross cap and 
$A\in I^{m-\frac14}(X,Y;C)$, then \linebreak$A^*A\in I^{2m-\frac12,\frac12}\left(\Delta_{T^*Y},\tilde{C}\right)$,where
$\tilde{C}\subset T^*Y\times  T^*Y$  intersects $\Delta_{T^*Y}$ cleanly in \linebreak codimension 1, and
$\tilde{C}$ is a folding canonical relation, i.e., both $\pi_L$ and $\pi_R$ are Whitney folds.
\end{thm}

It follows that the wavefront relation of $N:=A^*A$ is contained in $\Delta_{T^*Y}\cup \tilde{C}$,
with $\tilde{C}$ a folding canonical relation. By  \eqref{eqn HorSat}, \eqref{eqn iplwf} and that $\Delta_{T^*Y}$ 
acts as the identity relation on $T^*Y$, for any $u\in\mathcal E'(Y)$, we have
\be\label{eqn ccwf}
WF(Nu)\subseteq WF(u)\,  \cup\,  \left(\tilde{C}\circ WF(u)\right).
\ee
Furthermore, by \eqref{eqn ipl orders}, microlocally away from $\Delta_{T^*Y}\,\cap\, \tilde{C}$, we have
$N\in  I^{2m}\left(\Delta_{T^*Y}\setminus \tilde{C}\right)$, and
$N\in I^{2m-\frac12}\left(\tilde{C}\setminus \Delta_{T^*Y}\right)$; thus,
away from $\Delta_{T^*Y}\,\cap\, \tilde{C}$,  the order of the non-pseudodifferential operator
part of $N$, which constitutes an artifact, is $1/2$ lower order than the pseudodifferential  part of $N$.
Although the artifact's order  is 1/2 lower, the two-sided fold degeneracy of $\tilde{C}$,
combined with the paired Lagrangian nature of the normal operator,
produce a situation where it is not known whether the artifact is completely removable;
see \cite{FGP12} for further analysis and discussion. 

Under the assumptions of Theorem \ref{thm folded cc} and away from a very small microlocal set, 
this composition result and its implications 
apply to the dense array borehole data set, resulting in  artifacts 1/2 order smoother than the primary image.

%%%%%%%%%%%%%%%%%%%%%%%%%%%%%
%%%%%%%%%%%%%%%%% section break

\subsection{Proof of Thm. \ref{thm folded cc}}\label{subsec pf of folded cc}

Let $r_{min}<r_0<r_{max}$ and $\gamma_0\in \chi\left(\Lambda_{r_0}\right)\subset\Gamma^{c_0}$ 
be a Whitney fold point for $\pi_X:\Lambda_{r_0}\to \R^3_+$.
Repeating the analysis from Sec. \ref{subsec single}, we can assume that 
$x_1,x_2,p_3:=\xi_3|_{\Gamma_R^{c_0}}$ have
independent gradients near $\gamma_0$.
Since $\partial_r$ is transverse to $T\Lambda_r$,
$d\chi(\partial_r)$ is transverse to $T_{\gamma_0}\chi(\Lambda_r)$. Thus,
$(x_1,x_2,r,p_3)$ form  coordinates on $\Gamma^{c_0}$ near $\gamma_0$,
the acoustical distance function
described above  is a smooth function of  $(x_1,x_2,r,p_3)$, and
we can express $x_3$ and $(p_1, p_2):=(\xi_1,\xi_2)$ on $\Gamma^{c_0}$ in terms of them:
$x_3=f(x_1, x_2, r, p_3)$ and  $(p_1,p_2)=(g_1(x_1, x_2,r, p_3),
g_2(x_1, x_2, r, p_3))$ on $\Gamma^{c_0}$.
With respect to these coordinates,
$\pi_X(x_1, x_2, r, p_3)=(x_1, x_2, f(x_1, x_2, r, p_3))$  and \\

  $${d \pi_X}= \left(\begin{array}{cccc}
1 & 0 & 0 & 0 \\
0 & 1 & 0 & 0 \\
\frac{\partial f} {\partial x_1} & \frac{\partial f}{\partial x_2}
& \frac{\partial f}{\partial r} & \frac{\partial f}{\partial
p_3}
\end{array} \right).$$
\vskip .1 cm
\noindent From this we see that
$$\hbox{rank } d \pi_X = \left\{
\begin{array}{cc}
2, & {\rm if}  \ \ \frac{\partial f}{\partial r}=\frac{\partial f}{\partial p_3}=0\\
3,  & {\rm if}  \ \ \frac{\partial f}{\partial r} \neq 0   \
{\rm or} \ \frac{\partial f}{\partial p_3} \neq 0,
\end{array}  \right.$$
and $ \Sigma(\pi_X)=\{  \frac{\partial f}{\partial p_3}=
\frac{\partial f}{\partial r}=0 \}$. At points of $ \Sigma(\pi_X)$,  $\Ker d\pi_X$ is spanned by
$\{ (0,0, \delta r, \delta p_3): \delta r, \delta p_3\in\R  \}$, and the tangent space to $\Sigma(\pi_X)$
is
$$T \Sigma(\pi_X)= \Ker \left(d_{x_1,x_2,p_3,r}\left(\frac{\partial
f}{\partial p_3}\right)\right) \cap\Ker \left(d_{x_1,x_2,p_3,r}\left(\frac{\partial
f}{\partial r}\right)\right).$$
We have:
$$d_{x_1,x_2,p_3,r}\left(\frac{\partial
f}{\partial p_3}\right)= \left(\frac{\partial^2 f}{\partial x_1 \partial
p_3}, \frac{\partial^2 f}{\partial x_2
\partial p_3}, \frac{\partial^2 f}{
\partial p_3^2}, \frac{\partial^2 f}{\partial r \partial p_3}\right)$$ and
$$d_{x_1,x_2,p_3,r}\left(\frac{\partial
f}{\partial s}\right)= \left(\frac{\partial^2 f}{\partial x_1 \partial
r}, \frac{\partial^2 f}{\partial x_2
\partial r}, \frac{\partial^2 f}{
\partial p_3 \partial r}, \frac{\partial^2 f}{\partial r ^2}\right).$$
The assumption Def. \ref{def fold caustic}(ii) that $\pi_X: \Gamma^{c_0} \rightarrow \R^3_+$ 
is a submersion with folds implies that
 $\Sigma(\pi_X)$ is smooth
(i.e., these gradients are linearly independent), and
\begin{equation}\label{eqn ii}
T\Sigma(\pi_X)\hbox{  is
transverse to  }\Ker d \pi_X.
\end{equation}

We can parametrize
the canonical relation $C$ in terms of $r, x_1, x_2,
p_3$;  $(\alpha_1,\alpha_2)$,
where $(\alpha_1, \alpha_2, \sqrt {1-|\alpha|^2})$ is the unit the take off
direction of the reflected ray; and $\tau$, the variable dual  to time.
The incident ray travel time $t_{inc}$,  
the time it takes for a ray to travel from a source $S=(s,0)$ to an incident point $x$,
i.e., the acoustical distance from $S$ to $x$, 
can, by the nongrazing assumption  at the surface described in Sec. \ref{subsec model} and symmetry, be expressed in terms of $x$ and $\alpha$, $t_{inc}=t_{inc}(x,\alpha)$. 
On the other hand, the reflected ray travel time $t_{ref}$, 
the time it then takes for the reflected ray to reach the borehole at point $R=(0,0,r)$,  
is, again by symmetry, the acoustical distance from $R$ to $x$,  and thus 
$t_{ref}=t_{ref}(x_1,x_2,r,p_3)$;
the total  time for the single-reflection event is $t=t_{inc}+t_{ref}$.
Letting $(\rho, \sigma_1,\sigma_2,\tau)$ be the coordinates dual to $r, s_1,s_2,t$ in $T^*\D$,
we can take $(x_1,x_2,r,p_3,\alpha_1,\alpha_2,\tau)$ as coordinates on $C$, and
\bast
C &=& \Big\{ \big(r, s_1(x_1,x_2, f(x_1,x_2, r, p_3),
\alpha), s_2(x_1,x_2, f(x_1,x_2, r, p_3), \alpha),  \\
& & \quad t_{ref}(x_1, x_2,r,p_3) 
+ t_{inc}(x_1,x_2, f(x_1,x_2, r, p_3),\alpha),\\
& & \qquad \rho(x_1,x_2,f(x_1,x_2, r, p_3), p_3,\tau),
\sigma_1(x_1,x_2, f(x_1,x_2, r, p_3), \alpha,\tau),
\sigma_2(x_1,x_2,f(x_1,x_2, r, p_3), \alpha,\tau), \tau; \\
& &   \quad x_1, x_2, f(x_1, x_2,
r, p_3); -\tau\left(c_0^{-1}(x_1, x_2, f(x_1, x_2, r, p_3))
\alpha_1 + g_1(x_1, x_2, r, p_3)\right), \\
& &\qquad  -\tau\left(c_0^{-1}(x_1,
x_2, f(x_1, x_2, r, p_3)) \alpha_2 + g_2(x_1, x_2, r,
p_3)\right), \\
& & \qquad  -\tau\left(c_0^{-1}(x_1, x_2, f(x_1, x_2, r, p_3)\right)
\sqrt{1-|\alpha|^2} + p_3\big) \Big\},
\east
where $\rho(\cdot)$ is homogeneous of degree 1 in $p_3,\tau$, and
$\sigma_1(\cdot),\sigma_2(\cdot)$ are homogeneous of degree 1 in $\tau$.
\smallskip

\par We now show that $C$ is a folded cross cap in the sense of Def. \ref{def folded crosscap},
except possibly on a  set of codimension four.
In terms of the above coordinates on $C$ and the standard canonical coordinates $(x,\xi)$ on $T^*\R^3$, we 
can write $\pi_R$ as  a map $\pi_R : \R^7 \rightarrow \R^6$, given by
\bast
\pi_R(x_1,x_2,r,p_3,\alpha_1,\alpha_2,\tau)&=& \big ( x_1, x_2, f(x_1,x_2, r, p_3),  \\
& &  -\tau\left(c_0^{-1}\left(x_1, x_2, f\left(x_1, x_2, r,
p_3\right)\right) \alpha_1 + g_1(x_1, x_2, r, p_3)\right), \\
& & -\tau\left(c_0^{-1}\left(x_1, x_2, f\left(x_1, x_2, r, p_3\right)\right) \alpha_2 +
g_2(x_1, x_2, r, p_3)\right), \\
& &  -\tau\left(c_0^{-1}\left(x_1, x_2, f\left(x_1,x_2, r, p_3\right)\right) \sqrt{1-|\alpha|^2} + p_3\right) \big ).
\east

Thus,

 \be\label{matrix dpiR}
 {d \pi_R}= \left(\begin{array}{ccccccc}
1 & 0 & 0 & 0 & 0 & 0 & 0\\
0 & 1 & 0 & 0 & 0 & 0 & 0\\
\frac{\partial f}{\partial x_1} & \frac{\partial f}{\partial x_2} & \frac{\partial f}{\partial r} & \frac{\partial f}{\partial p_3} & 0 & 0 & 0\\
A_1 & A_2 & A_3 & A_4 & -\tau c_0^{-1} & 0 & -(c_0^{-1}\alpha_1+g_1)\\
B_1 & B_2 & B_3 & B_4 & 0 & -\tau c_0^{-1} & -(c_0^{-1}\alpha_2+g_2)\\
C_1 & C_2 & C_3 & C_4 & \tau c_0^{-1}
\frac{\alpha_1}{\sqrt{1-|\alpha|^2}} & \tau c_0^{-1}
\frac{\alpha_2}{\sqrt{1-|\alpha|^2} } & -(c_0^{-1} \sqrt{1-|\alpha|^2}+p_3)
\end{array} \right)
\ee
\noindent for some $A_j, B_j, C_j$.
The lower right $3\times 3$ submatrix
\vskip .5 cm
$$ \left(\begin{array}{ccc}
 -\tau c_0^{-1} & 0 & -(c_0^{-1}\alpha_1+g_1)\\
 0 & -\tau c_0^{-1} & -(c_0^{-1}\alpha_2+g_2)\\
 \tau c_0^{-1} \frac{\alpha_1}{\sqrt{1-|\alpha|^2}} & \tau c_0^{-1}
\frac{\alpha_2}{\sqrt{1-|\alpha|^2} } & -(c_0^{-1} \sqrt{1-|\alpha|^2}+p_3)
\end{array} \right)$$
is nonsingular, since its determinant
$$-\tau^2c_0^{-2}(1-|\alpha|)^{-\frac12}\left[ c_0^{-1}+\alpha_1 g_1
+\alpha_2 g_2+p_3(1-|\alpha|)^{\frac12}\right]$$
is nonzero: $(g_1,g_2,p_3)=\xi|_{\Gamma^{c_0}}$ and by  Assumption \ref{nsop},
any scattering over $\pi$ has been filtered out, i.e., 
$$\left(\alpha_1, \alpha_2 , (1-|\alpha|)^{\frac12}\right)\cdot
(p_1,p_2,p_3)\ne-c_0^{-1}.$$
It follows that

\be\label{eqn rank dpiR}
\hbox{ rank }  d \pi_R = \Big\{ \begin{array}{cc}
5, & {\rm if}  \ \ \frac{\partial f}{\partial p_3}= \frac{\partial f}{\partial r}=0,\\
6,  & {\rm if }  \ \ \frac{\partial f}{\partial p_3} \neq 0   \
{\rm or} \ \frac{\partial f}{\partial r} \neq 0.
\end{array} 
\ee
\medskip

\par Now,  $\Ker d \pi_R = \{ (0, 0, \delta r,  \delta p_3, \delta \alpha_1,
\delta \alpha_2, \delta \tau) \}$, where $\delta \alpha_1, \delta
\alpha_2, \delta \tau$ depend on $\delta p_3, \delta r$.
On the other hand, the tangent space to $\Sigma(\pi_R)$ is
$$T \Sigma(\pi_R)=\Ker
\left(d_{x_1,x_2,r,p_3,\alpha_1,\alpha_2,\tau}\left(\frac{\partial
f}{\partial p_3}\right)\right) \bigcap \Ker
\left(d_{x_1,x_2,r,p_3,\alpha_1,\alpha_2,\tau}\left(\frac{\partial
f}{\partial r}\right)\right),$$
where
$$d_{x_1,x_2,r,p_3,\alpha_1,\alpha_2,\tau}\left(\frac{\partial f}
{\partial p_3}\right)= 
\left(\frac{\partial^2 f}{\partial x_1 \partial p_3}, 
\frac{\partial^2 f}{\partial x_2 \partial p_3}, 
 \frac{\partial^2 f}{\partial r \partial p_3}, \frac{\partial^2 f}{\partial p_3^2},
0, 0,0\right)$$
and
$$d_{x_1,x_2,r,p_3,\alpha_1,\alpha_2,\tau}\left(\frac{\partial
f}{\partial r}\right)= 
\left(\frac{\partial^2 f}{\partial x_1 \partial r}, 
\frac{\partial^2 f}{\partial x_2 \partial s},
\frac{\partial^2 f}{\partial r ^2}, \frac{\partial^2 f}{\partial p_3 \partial r},  
 0, 0,0\right).$$
Combining  this with  \eqref{eqn ii},
one  sees that $\Ker d \pi_R$ is transverse to $T \Sigma(\pi_R)$
and thus, $\pi_R$ is a submersion with folds. 
One can also check that, off an exceptional set,  the image of the critical set is nonradial, i.e., 
$\xi\cdot dx \ne 0$ on $\Sigma(\pi_R)$. Since $T\left(\Sigma(\pi_R)\right)$ is the span of the 
columns of the matrix in \eqref{matrix dpiR} representing $d\pi_R$, 
while $\tau^{-1} \xi$ consists of  the last three entries in 
its last column,
and $\frac{\partial f}{\partial p_3} = \frac{\partial f}{\partial r}=0$ at $\Sigma(\pi_R)$, we see that
$(\xi\cdot dx)(W)=0$ for all $W\in T\left(\Sigma(\pi_R)\right)$ if and only if
\bast
c_0^{-1}\alpha_j+g_j+c_0^{-1}f_{x_j}\sqrt{1-|\alpha|^2}+f_{x_j}p_3=0,\, j=1,2,
\east
which defines a codimension 2 submanifold in $\Sigma(\pi_R)$ and hence codimension 4 in $C$.
\ss

\par We next need to show that $\pi_L$ is a cross cap.  
As for any canonical relation, $\Sigma(\pi_L)=\Sigma(\pi_R)$. 
Similar to the analysis for $\pi_R$, the projection $\pi_L: C\to T^*\D$ 
can be treated as mapping $\R^7 \rightarrow \R^8$, with (reordering the variables for convenience)
\bast
\pi_L(r, x_1, x_2, \alpha_1,\alpha_2, p_3, \tau)&\!\!\!=\!\!\!& \big(r, s_1(x_1,x_2, f(x_1,x_2, r, p_3),
\alpha),   s_2(x_1,x_2, f(x_1,x_2, r, p_3), \alpha), \\
& &  \quad t_{ref}(x_1, x_2,r,p_3) 
+ t_{inc}(x_1,x_2, f(x_1,x_2, r, p_3),\alpha), \\
& & \quad \rho\left(x_1,x_2, f(x_1, x_2, r, p_3), p_3, \tau\right),  \\
& & \quad
\sigma_1(x_1,x_2, f(x_1,x_2, r, p_3), \alpha, \tau), \sigma_2(x_1,x_2,
f(x_1,x_2, r, p_3), \alpha, \tau), \tau \Big )
\east
and thus,
\smallskip

${d \pi_L}= \left(\begin{array}{ccccccc}
1 & 0 & 0 & 0 & 0 & 0 & 0\\
\frac{\partial s_1}{\partial x_3}\frac{\partial f}{\partial r} & \frac{\partial s_1}{\partial
x_1} +  \frac{\partial s_1}{\partial x_3}\frac{\partial
f}{\partial x_1} & \frac{\partial s_1}{\partial x_2} +
\frac{\partial s_1}{\partial x_3}\frac{\partial f}{\partial x_2} &
 \frac{\partial s_1}{\partial \alpha_1} & \frac{\partial s_1}{\partial \alpha_2} & \frac{\partial s_1}{\partial x_3} \frac{\partial f}{\partial p_3} & 0\\
\frac{\partial s_2}{\partial x_3}\frac{\partial f}{\partial r} & \frac{\partial s_2}{\partial
x_1} + \frac{\partial s_2}{\partial x_3}\frac{\partial f}{\partial
x_1} & \frac{\partial s_2}{\partial x_2}  +  \frac{\partial
s_2}{\partial x_3}\frac{\partial f}{\partial x_2}& \frac{\partial
s_2}{\partial \alpha_1}
 & \frac{\partial s_2}{\partial \alpha_2} & \frac{\partial s_2}{\partial x_3} \frac{\partial f}{\partial p_3} & 0\\
\frac{\partial t_{ref}}{\partial r}
+\frac{\partial t_{inc}}{\partial x_3}\frac{\partial f}{\partial r} & \frac{\partial
t_{ref}}{\partial x_1} + \frac{\partial t_{inc}}{\partial
x_3}\frac{\partial f}{\partial x_1} & \frac{\partial
t_{ref}}{\partial x_2} +  \frac{\partial t_{inc}}{\partial
x_3}\frac{\partial f}{\partial x_2} & \frac{\partial
t_{inc}}{\partial \alpha_1} &
 \frac{\partial t_{inc}}{\partial \alpha_2} & \frac{\partial t_{ref}}{\partial p_3} + \frac{\partial t_{inc}}{\partial x_3}
\frac{\partial f}{\partial p_3} & 0 \\
\frac{\partial \rho}{\partial r} & 
\frac{\partial \rho}{\partial x_1} +  \frac{\partial \rho}{\partial x_3}\frac{\partial f}{\partial x_1} & 
\frac{\partial\rho}{\partial x_2} + \frac{\partial \rho}{\partial x_3}\frac{\partial f}{\partial x_2}& 0 & 0 
& \frac{\partial \rho}{\partial x_3} \frac{\partial f}{\partial p_3} + \frac{\partial \rho}{\partial p_3}& \frac{\partial \rho}{\partial \tau} \\
\frac{\partial \sigma_1}{\partial x_3}\frac{\partial f}{\partial r} & \frac{\partial
\sigma_1}{\partial x_1} +  \frac{\partial \sigma_1}{\partial
x_3}\frac{\partial f}{\partial x_1}& \frac{\partial
\sigma_1}{\partial x_2}  +  \frac{\partial \sigma_1}{\partial
x_3}\frac{\partial f}{\partial x_2} & \frac{\partial
\sigma_1}{\partial \alpha_1} & \frac{\partial \sigma_1}{\partial
\alpha_2} &
\frac{\partial \sigma_1}{\partial x_3} \frac{\partial f}{\partial p_3} & \frac{\partial \sigma_1}{\partial \tau} \\
\frac{\partial \sigma_2}{\partial x_3}\frac{\partial f}{\partial r} & \frac{\partial
\sigma_2}{\partial x_1} +  \frac{\partial \sigma_2}{\partial
x_3}\frac{\partial f}{\partial x_1}& \frac{\partial
\sigma_2}{\partial x_2} +  \frac{\partial \sigma_2}{\partial
x_3}\frac{\partial f}{\partial x_2} & \frac{\partial
\sigma_2}{\partial \alpha_1} & \frac{\partial \sigma_2}{\partial
\alpha_2} &
\frac{\partial \rho_2}{\partial x_3} \frac{\partial f}{\partial p_3} & \frac{\partial \sigma_2}{\partial \tau} \\
0 & 0 & 0 & 0 & 0 & 0 & 1
\end{array}  \right)$.
\smallskip

\par

Since corank $ d\pi_L=$ corank $d\pi_R$,  it follows from \eqref{eqn rank dpiR} that
$d\pi_L$ is injective except where it has a one-dimensional kernel above  caustic
points, where $\frac{\partial f}{\partial p_3}=\frac{\partial f}{\partial r}=0$.
We will use these two conditions, plus one more, in order to simplify this matrix.
Namely, by rotation about the borehole, we can assume that, for the point of interest,  the tangent plane
$\Pi=\pi_X\left(T_{\lambda_0}\Lambda_r^{c_0}\right)$ from Sec. \ref{subsec single} is the graph of $x_3$ as a linear function independent of $x_2$. As a consequence, $f_{x_2}=0$ at this point. The matrix for $d\pi_L$ then becomes 
\ss

\be\label{matrix dpiL simple}
{d \pi_L}= \left(\begin{array}{ccccccc}
1 & 0 & 0 & 0 & 0 & 0 & 0\\
0 &  \frac{\partial s_1}{\partial
x_1} +  \frac{\partial s_1}{\partial x_3}\frac{\partial
f}{\partial x_1} & \frac{\partial s_1}{\partial x_2} &
 \frac{\partial s_1}{\partial \alpha_1} & \frac{\partial s_1}{\partial \alpha_2} & 0 & 0\\
0 & \frac{\partial s_2}{\partial
x_1} + \frac{\partial s_2}{\partial x_3}\frac{\partial f}{\partial
x_1} & \frac{\partial s_2}{\partial x_2} & \frac{\partial
s_2}{\partial \alpha_1}
 & \frac{\partial s_2}{\partial \alpha_2} & 0 & 0\\
\frac{\partial t_{ref}}{\partial r}
 & \frac{\partial
t_{ref}}{\partial x_1} + \frac{\partial t_{inc}}{\partial
x_3}\frac{\partial f}{\partial x_1} & 
\frac{\partial t_{ref}}{\partial x_2}  & 
\frac{\partial t_{inc}}{\partial \alpha_1} &
 \frac{\partial t_{inc}}{\partial \alpha_2} & \frac{\partial t_{ref}}{\partial p_3}  & 0 \\
\frac{\partial \rho}{\partial r} & 
\frac{\partial \rho}{\partial x_1} +  \frac{\partial \rho}{\partial x_3}\frac{\partial f}{\partial x_1}   & 
\frac{\partial\rho}{\partial x_2} & 0 & 0 &  \frac{\partial \rho}{\partial p_3} & \frac{\partial \rho}{\partial \tau} \\
0 & \frac{\partial
\sigma_1}{\partial x_1} +  \frac{\partial \sigma_1}{\partial
x_3}\frac{\partial f}{\partial x_1}& \frac{\partial
\sigma_1}{\partial x_2}   & \frac{\partial
\sigma_1}{\partial \alpha_1} & \frac{\partial \sigma_1}{\partial
\alpha_2} &
0 & \frac{\partial \sigma_1}{\partial \tau} \\
0 & \frac{\partial
\sigma_2}{\partial x_1} +  \frac{\partial \sigma_2}{\partial
x_3}\frac{\partial f}{\partial x_1}& \frac{\partial
\sigma_2}{\partial x_2} & \frac{\partial
\sigma_2}{\partial \alpha_1} & \frac{\partial \sigma_2}{\partial
\alpha_2} &
0 & \frac{\partial \sigma_2}{\partial \tau} \\
0 & 0 & 0 & 0 & 0 & 0 & 1
\end{array}  \right).
\ee
\smallskip

Writing a spanning element $V_L\in \Ker d \pi_L$ as $(\delta r, \delta x_1,\delta x_2, \delta \alpha_1,
\delta \alpha_2, \delta p_3,  \delta \tau)$, we see from \eqref{matrix dpiL simple} that
$\delta r=\delta\tau=0$, 
so that the $\frac{\partial \sigma_1}{\partial \tau},\, \frac{\partial \sigma_2}{\partial \tau}$ terms in 
\eqref{matrix dpiL simple} can be ignored. 

Furthermore, under the assumption that the rays from the sources are transverse to the caustic surface 
(see \cite{N00}),  the determinant

\be\label{eqn minor}
\left|\begin{array}{cccc}
   \frac{\partial s_1}{\partial x_1}&   \frac{\partial s_1}{\partial x_2}  &  \frac{\partial s_1}{\partial \alpha_1} &  \frac{\partial s_1}{\partial \alpha_2} \\
 \frac{\partial s_2}{\partial x_1}&   \frac{\partial s_2}{\partial x_2}  &  \frac{\partial s_2}{\partial \alpha_1} &  \frac{\partial s_2}{\partial \alpha_2} \\
 \frac{\partial \sigma_1}{\partial x_1}&   \frac{\partial \sigma_1}{\partial x_2}  &  \frac{\partial \sigma_1}{\partial \alpha_1} &  \frac{\partial \sigma_1}{\partial
\alpha_2}\\
\frac{\partial \sigma_2}{\partial x_1}& \frac{\partial
\sigma_2}{\partial x_2}  &  \frac{\partial \sigma_2}{\partial
\alpha_1} & \frac{\partial \sigma_2}{\partial \alpha_2}
\end{array} \right|\ne 0.
\ee 
The matrix in \eqref{eqn minor} is almost a minor of 
\eqref{matrix dpiL simple}: they differ only in the first column, by
\be\label{eqn matrix diff}
f_{x_1}\cdot\left[(s_1)_{x_3}, (s_2)_{x_3}, (\sigma_1)_{x_3}, (\sigma_2)_{x_3}\right]^T.
\ee
Thus, if we make the small slope assumption that 
\be\label{cond dip}
|f_{x_1}|\quad\hbox{  is sufficiently small,}
\ee
i.e., the normal to the fold caustic surface is sufficiently close to vertical,
then the corresponding minor of \eqref{matrix dpiL simple} is nonsingular,
which implies that
$\delta x_1=\delta x_2=\delta \alpha_1=\delta \alpha_2=0$.
Thus, $\delta p_3$ is the only nonzero entry in $V_L$; note also that this then implies that 
$\frac{\partial t_{ref}}{\partial p_3}= \frac{\partial \rho}{\partial p_3}=0$ at caustics. 
It then  follows from \eqref{eqn causticfold conds}  that (iii) below Def. \ref{def crosscap} is satisfied, and (i), (ii) follow from this analysis as well. Hence, $\pi_L$
is a cross cap. 
Again, one can check that, away from a set of high codimension,  
the image of the cross cap points is nonradial in 
$T^*\mathbb D$.
A point is radial if and only if 
$\left(\rho,\sigma_1,\sigma_2,\tau\right)^T\cdot W=0$
for all $W\in T\left(\pi_L(\Sigma(\pi_L))\right)$, 
i.e., for $W$ 
in the span of the columns of the upper $4\times 7$ submatrix of $d\pi_L$,
this becomes

\bast
\rho &=& -(t_{ref})_r\tau\\
((s_1)_{x_1}+\epsilon_1)\sigma_1+(s_2)_{x_1}\sigma_2 &=& -(t_{ref})_{x_1}\tau\\
((s_1)_{x_2}+\epsilon_2)\sigma_1+(s_2)_{x_2}\sigma_2 &=& -(t_{ref})_{x_2}\tau\\
((s_1)_{\alpha_1}+\epsilon_3)\sigma_1+ (s_2)_{\alpha_1}\sigma_2 &=& -(t_{inc})_{\alpha_1}\tau\\
((s_1)_{\alpha_2}+\epsilon_4)\sigma_1+ (s_2)_{\alpha_2}\sigma_2 &=& -(t_{inc})_{\alpha_2}\tau,
\east
for some small $\epsilon_j,\, 1\le j\le 4$.
The first equation imposes one condition.
On the other hand, the last four equations impose two more,
since the coefficient matrix is the upper $2\times 4$ submatrix of the matrix in \eqref{eqn minor} and thus
has rank two, meaning that the right hand sides of these last four equations must satisfy two linear conditions
in order for the equations to be solvable. Thus, the set of possibly radial points of
$\pi_L\left(\Sigma(\pi_L)\right)$ is of codimension at least  three in the critical set, and thus codimension 5 in 
$C$. Combined with the codimension 4 set of possible radial points of $\pi_R$, we see that the nonradiality 
conditions of Def. \ref{def folded crosscap} are satisfied away from a set of codimension at least 4 in $C$.
\smallskip

In summary, we have shown that if the  ray geometry of the background sound speed $c_0$
 has at most fold caustics with
respect to the borehole, and the small slope assumption \eqref{cond dip} holds,
 then away from a codimension 4 set the canonical relation $C$ is a folded cross cap, 
 finishing the proof of  Thm. \ref{thm folded cc}. \hspace{5cm}$\Box$
  \smallskip
 
Theorem \ref{thm folded cc} then implies that Theorem \ref{thm folded crosscap} applies to the composition 
forming the normal operator $F^*F$
(away from the possible bad set microlocally),
with the consequences for artifacts as described above.

%%%%%%%%%%%%%%%%%%%%%%%%%%%%%
%%%%%%%%%%%%%%%%% section break

\section{Crosswell and walkaway geometries}\label{crosswell}

As for the dense array in Sec. \ref{sec dense constant}, for the crosswell and walkaway geometries we   
compute  $F=d\F$
at the constant background sound speed $c_0=1$
 by restricting the basic phase function
\eqref{eqn phi basic} to each data set $\D$.

%%%%%%%%%%%%%%%
\subsection{Crosswell geometry}\label{subsec crosswell}

For the crosswell  (CW) geometry, we assume that the sources and receivers are located in
parallel, vertical boreholes. For simplicity, assume that the  sources form an open interval along
 the line $y_1=s_0$, $y_2 =0$,  for some $s_0>0$,

\[\Sigma_S=\{(s_0, 0, s):\, s\in (s_{min},s_{max})=:I_S\},\]
and the  receivers are similarly located, as for the other geometries, in the borehole along the $y_3$ axis, say
\[\Sigma_R=\{(0,0,r):\, r\in (r_{min},r_{max})=:I_R\},\]
and we  identify  $\D=\D_{CW}=\left(s_{min},s_{max}\right)\times \left(r_{min},r_{max}\right)
\times (t_{min},t_{max})$.
\smallskip

The associated linearized  scattering operator $F$ is then a Fourier integral operator with
phase function
obtained by restricting \eqref{eqn phi basic} to $\D_{CW}$:

\begin{equation}\label{phase cross-well}\nonumber
\phi_{CW}(s,r,t,y;\omega) = \left(t-\sqrt{(y_1-s_0)^2+y_2^2+(y_3-s)^2}
-\sqrt{(y_1^2+y_2^2+(y_3-r)^2} \right) \omega.
\end{equation}

The structure of  the linearized scattering operator $F$ for the crosswell geometry is summarized by the following.
\ms

\begin{theorem}\label{thm crosswell}
The linearized scattering operator $F$ for the crosswell imaging
\lb geometry is a  Fourier integral operator, $F\in I^{-\frac12}\left(C_{CW}\right)$,
whose canonical relation $C_{CW}$
is singular on
the union of two hypersurfaces, $\Sigma^1 \cup \Sigma^2$, with
$\Sigma^1$ and $\Sigma^2$ intersecting transversally.
On  $\Sigma^1\setminus\Sigma^2$,
$\pi_L$ has a fold  singularity and $\pi_R$ is a blowdown, while on
$\Sigma^2\setminus\Sigma^1$, both of the projections $\pi_L$ and $\pi_R$
have fold singularities.
\end{theorem}

\begin{proof}[Proof of Thm. 6.1]
Let
$$ A :=\sqrt{(y_1 - s_0)^2+y_2^2+(y_3-s)^2} \, , \quad
B := \sqrt{y_1^2+y_2^2+(y_3-r)^2}.$$

We calculate the canonical relation, $C_{CW}$, parametrized by $\phi_{CW}$, 
and  classify the singularities of the left and
right projections.  We have:

\begin{eqnarray}
C_{CW}&=&\bigg\{\Big(s, r, A+B, \frac{y_3-s}{A} \omega,  \frac{y_3-r}{B} \omega, \omega;\nonumber\\
& & \quad y_1, y_2, y_3, \left(\frac{y_1 - s_0}{A}+  \frac{y_1}{B}\right)\omega,
\left(\frac{y_2}{A}+  \frac{y_2}{B}\right)\omega, \left(\frac{y_3-s}{A}+
\frac{y_3-r}{B}\right)\omega\Big)\nonumber\\
& & \qquad : y\in\R^3,\, s\in I_S,\, r\in I_R,\, \omega\ne 0\bigg\}.\nonumber
\end{eqnarray}

\par With respect to these coordinates,  the left projection, $\pi_L:C_{CW}\to T^*\D_{CW}$, is

\begin{equation}\label{piL}\nonumber
\pi_L(y,s,r, \omega)=\Big(s, r, \omega, A+B,  \frac{y_3-s}{A} \omega,  \frac{y_3-r}{B} \omega\Big)
\end{equation}

and the right projection, $\pi_R:C_{CW}\to T^*\R^3$, is

\begin{equation}\label{piR}\nonumber
\pi_R(y, s, r, \omega)=\Big(y_1, y_2, y_3, \left(\frac{y_1 - s_0}{A}+  \frac{y_1}{B}\right)\omega,  \left(\frac{y_2}{A}+  \frac{y_2}{B}\right)\omega, \left(\frac{y_3-s}{A}+  \frac{y_3-r}{B}\right)\omega\Big).
\end{equation}

\par We first study $\pi_L$. Denote the variables dual to $s,r,t$ by $\sigma, \rho, \tau$, resp.
Since $\pi_L$ is the identity in the $s,r, \omega$ variables,  det $d\pi_L$ equals
det $D(\sigma,\rho,\tau)/Dy$,  i.e.,

\begin{eqnarray}\label{det piL}
\textnormal{det}(d\pi_L) &=& \left| \begin{array} {ccc}  \frac{y_1-s_0}{A}+  \frac{y_1}{B}& \frac{y_2}{A}
+  \frac{y_2}{B} & \frac{y_3-s}{A}+  \frac{y_3-r}{B} \\
-\frac{(y_3-s)(y_1-s_0)}{A^3}\omega & -\frac{(y_3-s)y_2}{A^3}\omega & \frac{(y_1 - s_0)^2+y_2^2}{A^3}\omega\\
-\frac{(y_3-r)y_1}{B^3}\omega & -\frac{(y_3-r)y_2}{B^3}\omega & \frac{y_1^2+y_2^2}{B^3}\omega\\
\end{array} \right|\nonumber\\
& & \nonumber \\
&=& -\frac{\omega^2}{A^3B^3} s_0 y_2 \Big(\frac{y_3-s}{A}+  \frac{y_3-r}{B}\Big).\nonumber
\end{eqnarray}
Thus, $\hbox{det}(d\pi_L)=0$  on $\Sigma^1 \cup \Sigma^2$, where

\[
\Sigma^1 := \left\{ y_2=0 \right\},\qquad
\Sigma^2 := \left\{ \frac{y_3-s}{A}+  \frac{y_3-r}{B}=0 \right\}.
\]
\smallskip

Note that points in $\Sigma^1 \cap \Sigma^2$ correspond to unbroken rays from $\Sigma_S$ to $\Sigma_R$,
not undergoing any scattering, and thus are first arrival events.
One can thus filter the data away from $\Sigma^1 \cap \Sigma^2$ by multiplying  $d(s,r,t)$
by a smooth cutoff $\chi\left(c_0t-\left|S(s)-R(r)\right|\right)$, 
where $supp\left(\chi\right)\subseteq \{t\ge \epsilon\}$
for some $\epsilon>0$. Hence, we do not need to consider the more singular structure
of $C_{CW}$ at $\Sigma^1 \cap \Sigma^2$.
\ms

Along each of $\Sigma^1\setminus \Sigma^2$ and $\Sigma^2\setminus\Sigma^1$,
$\det d\pi_L$ vanishes simply, and thus $d\pi_L$ drops rank
by 1. One easily  sees that,  along $\Sigma^1\setminus \Sigma^2$,
 $\Ker d\pi_L=\frac{\partial}{\partial y_2}$
 and hence (cf. Def. \ref{def fold}) $\pi_L$ has a {fold singularity}
 at points of $\Sigma^1\setminus\Sigma^2$.
 Similarly, $\Ker d\pi_L=\frac{\partial}{\partial y_2}$ at points of $\Sigma^2\setminus\Sigma^1$,
 and hence $\pi_L$ has a fold singularity there as well.
 \ms

\par Next, we consider $\pi_R$. As for any canonical relation, $d\pi_R$ also drops rank by
the same amount as $d\pi_L$,  and hence by $1$
on $\left(\Sigma^1\setminus\Sigma^2\right)\cup \left(\Sigma^2\setminus\Sigma^1\right)$.
We find its kernel by computing

\begin{eqnarray}\label{det piR}
\frac{D\eta}{D(s,r,\omega)} &=& \left[ \begin{array} {ccc}  \frac{(y_1 -s_0)(y_3 -s)}{A^3}\:\omega & \frac{y_1(y_3 -r)}{B^3}\:\omega  & \frac{y_1 - s_0}{A}+  \frac{y_1}{B} \\
\frac{y_2(y_3-s)}{A^3}\omega & \frac{y_2(y_3-r)}{B^3}\omega & \frac{y_2}{A} + \frac{y_2}{B}\\
-\frac{(y_1-s_0)^2 + y_2^2}{A^3}\omega & -\frac{y_1^2 + y_2^2}{B^3}\omega & \frac{y_3 -s}{A} + \frac{y_3 -r}{B}\\
\end{array} \right].
\end{eqnarray}

The kernel of $d\pi_R$ is contained in
$span\left\{\frac{\partial}{\partial s}, \frac{\partial}{\partial r},
\frac{\partial}{\partial \omega}\right\}$, which when applied to the defining function
$y_2$ of $\Sigma^1$
gives $0$. Hence, along $\left(\Sigma^1\setminus\Sigma^2\right)$,
ker($d\pi_R)\subset T\Sigma^1$ and thus (cf. Def. \ref{def blowdown})
$\pi_R$ has  a {blowdown singularity} along
$\Sigma^1\setminus \Sigma^2$.
On the other hand,  along $\Sigma^2\setminus\Sigma^1$,  the first 2 entries of the last row
are nonzero while the last one is $0$. Hence the kernel of $d\pi_R$ is spanned by
$\frac{\partial}{\partial s}$ or  $\frac{\partial}{\partial r}$, which is transverse to $\Sigma^2$,
and so $\pi_R$ has a {fold singularity}.

\end{proof}

%%%%%%%%%%%%%%%%%%%%%%%%%%%%%
%%%%%%%%%%%%%%%%% section break

\subsection{Walkaway geometry}\label{subsec walkaway}

For the walkaway geometry,  the set of sources is assumed to be an open subset of the $y_1$ axis,
\[\Sigma_S=\{(s,0,0):\, s\in I_S=(s_{min},s_{max})\}\]
and the set of receivers is as throughout an open subset of the $y_3$ axis,
\[\Sigma_R=\{(0,0,r):\, r\in I_R=(r_{min},r_{max})\},\]
so that  $\D=\D_{WA}=I_S\times I_R\times I_T$.
Restricting \eqref{eqn phi basic} to $\D_{WA}$,
the phase function of $F$  is
\[\phi_{WA}(s,r,t,y;\omega) = \left(t- \sqrt{(y_1 -s)^2 + y_2^2 + y_3^2} - \sqrt{y_1^2 + y_2 ^2 
+ (y_3 - r)^2}\right)\omega.\]
Let
\begin{eqnarray}
A &:=& \sqrt{(y_1-s)^2+y_2^2+y_3^2};\label{A-WA}\\
B &:=& \sqrt{y_1^2+y_2^2+(y_3 -r)^2}\label{B-WA}.
\end{eqnarray}

The structure of  the linearized scattering operator $F$ for the walkaway geometry is summarized by
the following.
\ms

\begin{theorem}\label{thm walkaway}
The linearized scattering operator $F$ for the walkaway geometry is
a Fourier integral operator, $F\in I^{-\frac12}\left(C_{WA}\right)$,
 whose canonical relation $C_{WA}$
is singular at the union of two smooth hypersurfaces $\Sigma^1$ and $\Sigma^2$, which intersect transversally.
At $\Sigma^1\setminus \Sigma^2$,
$\pi_L$ has a fold singularity  and $\pi_R$ has a blowdown singularity,
while at $\Sigma^2\setminus \Sigma^1$, $\pi_L$ is a fold at all points and $\pi_R$ is a fold
away from a hypersurface.
\end{theorem}
\ms

\begin{proof}[Proof of Thm. 6.2] 
The canonical relation $C_{WA}$ of $F$  is

\begin{eqnarray}
C_{WA}&=&\bigg\{\Big(s, r, A+B, \frac{y_1-s}{A} \omega,  \frac{y_3 -r}{B} \omega, \omega;\nonumber\\
& & \quad y_1, y_2, y_3, \left(\frac{y_1 -s }{A}+\frac{y_1}{B}\right)\omega,
\left(\frac{y_2}{A}+  \frac{y_2}{B}\right)\omega, \left(\frac{y_3}{A}+\frac{y_3-r}{B}\right)\omega  \Big)\nonumber\\
& & \qquad : y\in\R^3,\, s\in I_S,\, r\in I_R,\, \omega\ne 0\bigg\}.\nonumber
\end{eqnarray}

The right projection $\pi_R:C_{WA}\to T^*\R^3$ is
\begin{equation}\label{piR}
\pi_R(y,s, r, \omega)=\bigg(y_1, y_2, y_3, \left( \frac{y_1 - s}{A}+ \frac{y_1}{B}\right)\omega,  \left(\frac{y_2}{A}
+  \frac{y_2}{B}\right)\omega, \left(\frac{y_3}{A}+ \frac{y_3-r}{B} \right)\omega\bigg).
\end{equation}

Since $\pi_R$ is the identity in the $y$ variables, to  compute the det $d\pi_R$ we only need to
compute the Jacobian in the remaining variables $s, r, \omega$, which (in this order) is

\begin{eqnarray}\label{eqn dpiR}
\frac{D(\eta_1,\eta_2,\eta_3)}{D(s,r,\omega)}&=& \left[ \begin{array} {ccc}
-\frac{y_2^2 + y_3^2}{A^3}\omega &  \frac{y_1(y_3 -r)}{B^3}\omega &  \frac{y_1-s}{A} + \frac{y_1}{B}\\
\frac{y_2(y_1-s)}{A^3}\omega & \frac{y_2(y_3 -r) }{B^3}\omega & \frac{y_2}{A}+\frac{y_2}{B}\\
\frac{y_3(y_1-s)}{A^3}\omega & -\frac{y_1^2 + y_2^2}{B^3}\omega & \frac{y_3}{A}+\frac{y_3 -r}{B} \\
\end{array} \right].\\\nonumber
\end{eqnarray}
A calculation yields that
\begin{eqnarray}
\textnormal{det}\left(d\pi_R\right)
 &=& {{-}}\frac{\omega^2 y_2}{A^2 B^2} \left(\frac{y_1 (y_1 -s) + y_2^2 + y_3^2}{A} 
 + \frac{y_1^2 + y_2^2 + y_3 (y_3 -r)}{B}\right);\nonumber
\end{eqnarray}
the expression in the parentheses  can be written as $y\cdot\left(\frac{y-S}{A}+ \frac{y-R}{B}\right)$.
If we let
\be\label{def Sigmas}
\Sigma^1 := \left\{f_1:= y_2=0 \right\},\qquad
\Sigma^2 := \left\{ f_2:=y\cdot\left(\frac{y-S}{A}+ \frac{y-R}{B}\right)=0 \right\},
\ee
then   $\Sigma^1,\,  \Sigma^2$ intersect transversally.
Furthermore,  on $\left(\Sigma^1 \setminus \Sigma^2\right) \cup \left(\Sigma^2 \setminus \Sigma^1\right)$, 
$\det d\pi_R$ vanishes  simply,  
and thus $d\pi_R$  drops rank by 1 there;
by general principles concerning canonical relations, 
the same facts hold for $\det d\pi_L=\det d\pi_R$ and rank $d\pi_L$, resp.
\smallskip

\par  $\Sigma^1\setminus\Sigma^2$: From \eqref{piR} we see that  
$\Ker d\pi_R\subset span\left\{\frac{\partial}{\partial s}, 
\frac{\partial}{\partial r}, \frac{\partial}{\partial \omega}\right\}$, 
which is contained in $ T\Sigma^1$,
and  is one-dimensional at points of $\Sigma^1\setminus \Sigma^2$; 
hence, $\pi_R$ has  a {blowdown singularity} there.
Next  consider  $\pi_L:C_{WA}\to T^*\D_{WA}$,
\begin{equation}\label{piL}
\pi_L (y,s,r,\omega) = \left(s,r, A+B ,\frac{y_1 -s}{A}\:\omega, \frac{y_3 -r}{B}\:\omega,  \omega\right).
\end{equation}
As noted above, $d\pi_L$ drops rank by the same amount as $d\pi_R$ and so also has a one-dimensional
kernel  along $\Sigma^1\setminus \Sigma^2$. Since $\pi_L$ is  the identity in
the $s,r,\omega$ variables, we only need  compute the differential in the remaining variables, $y$,
\begin{equation}\label{dpiL}
\frac{D(t,\sigma,\rho)}{D(y_1,y_2,y_3)}
= \left[ \begin{array} {ccc}
\frac{(y_1-s)}{A}+\frac{y_1}{B} & \frac{y_2}{A}+\frac{y_2}{B} & \frac{y_3}{A}+ \frac{y_3 -r}{B}\\
\frac{y_2^2 + y_3^2}{A^3}\:\omega &  -\frac{y_2(y_1 -s)}{A^3}\:\omega & -\frac{y_3 (y_1-s)}{A^3}\:\omega \\
-\frac{y_1(y_3-r)}{B^3}\:\omega & -\frac{y_2(y_3 -r) }{B^3}\:\omega & \frac{y_1^2 + y_2^2}{B^3}\:\omega\\
\end{array} \right],
\end{equation}
and $\Ker d\pi_L$ is contained in span $\{ \frac{\partial}{\partial y_1}, 
\frac{\partial}{\partial y_2}, \frac{\partial}{\partial y_3} \}$.
Since the entries in the middle column are multiples of $y_2$, which vanishes on $\Sigma^1$, 
one sees that, on $\Sigma^1\setminus \Sigma^2$,
$\Ker d\pi_L=span\{\frac{\partial}{\partial y_2}\}$,
which is transverse to $\Sigma^1=\{y_2=0\}$. Thus, $\pi_L$ has a {fold singularity} 
along $\Sigma^1\setminus\Sigma^2$.
\medskip

$\Sigma^2\setminus \Sigma^1$: We show that all the singularities of $\pi_L$  are of fold type,
while $\pi_R$ has  fold singularities on the complement of a subset defined by a polynomial equation.
 At points of $\Sigma^2\setminus\Sigma^1$, as was the case on $\Sigma^1\setminus\Sigma^2$, 
$\Ker d\pi_L\subset\{ \frac{\partial}{\partial y_1}, \frac{\partial}{\partial y_2}, \frac{\partial}{\partial y_3} \}$ and  is 
one-dimensional. If  $V_L \ne 0$ spans $\Ker d\pi_L$, then it is annihilated by all of the rows of \eqref{dpiL}, 
 in particular the first row, and thus   $\la d_y(A+B), V_L\ra=0$.

\par There is a geometric interpretation of this last fact: for fixed $s,r$,  
the family of level surfaces of $A+B$, $E_{s,r,t}:=\{y: A+B=t \}$, indexed by $t > \sqrt{s^2+t^2}$,
are ellipsoids  with foci at $s$ and $r$, and outward pointing (nonunit) normal $\nu:=d_{y}(A+B)$. 
Then, since $\la \nu, V_L\ra=0$, we see that $V_L$ is tangent to $E_{s,r,t}$;
on the other hand, by \eqref{def Sigmas}, $y$ (considered as a vector) is also tangent to $E_{s,r,t}$.

\par Notice that, with $f_2=y\cdot \nu$  as in \eqref{def Sigmas},  $d_yf_2=\nu + y^t d_y \nu$. One has 
\begin{equation}\label{new eq}
\la d_yf_2, V_L \ra= \la\nu, V_L \ra  + y^t (d_y \nu) V_L.
\end{equation}
The first term on the right hand side of \eqref{new eq} is zero and the second one is positive 
since the ellipsoid $E_{s,r,t}$ has positive curvature (for every 
$V, V' \in T_y E_{s,r,t}, \ V^t (d_y\nu)V' > 0$) and $y,\, \nu\in T_yE_{s,r,t}$. 
Thus $\pi_L$ has a fold singularity along $\Sigma^2$.
\ss

For $\pi_R$, 
$\Ker d\pi_R\subset span\left\{\frac{\partial}{\partial s}, 
\frac{\partial}{\partial r}, \frac{\partial}{\partial \omega}\right\}$  is one-dimensional, and thus  spanned by a
$V_R=\delta s\frac{\partial}{\partial s}+ \delta r \frac{\partial}{\partial r}
+ \delta \omega \frac{\partial}{\partial \omega}$. 
From the matrix \eqref{eqn dpiR} representing the essential part of $d\pi_R$, we use the second row to solve for 
$\delta \omega$ in terms of $\delta s$  and $\delta  r$ (the value of which will be irrelevant below), and the first 
row to solve for $\delta s$ in terms of 
$\delta r$, namely $\delta s=-\delta r \frac{A}{B}  \frac{s}{r}$. 
Thus, $V_R=-\delta r \frac{A}{B}  \frac{s}{r} \frac{\partial}{\partial s}+ \delta r \frac{\partial}{\partial r}
+ \delta \omega \frac{\partial}{\partial \omega}$; applying this to $f_2$ (which is independent of $\omega$),
a calculation gives  the critical set 
\be\label{eqn Zariski}
\Sigma\left(\pi_R|_{\Sigma^2\setminus \Sigma^1}\right)
=\left\{ s^2B^2(y_2^2+y_3^2)-r^2 A^2 (y_1^2+y_2^2)=0 \right\}.
\ee 
We can see that the polynomial defining function in \eqref{eqn Zariski} is nonzero at some points, 
e.g., by taking $s$ or $r\to\infty$ and considering the leading coefficient in $s$ or $r$, resp.
Therefore, $\Sigma\left(\pi_R|_{\Sigma^2\setminus \Sigma^1}\right)$ is a lower dimensional variety, 
whose complement  in $\Sigma^2\setminus \Sigma^1$ is  dense;
on that set,  $V_Rf_2\ne0$ so that $\pi_R$  has a fold singularity at those points. 
This finishes the proof of Thm. \ref{thm walkaway}.
\end{proof}

%%%%%%%%%%%%%%%%%%%%%%%%%%%%%
%%%%%%%%%%%%%%%%% section break

\subsection{Artifacts for crosswell and walkaway}\label{subsec artifacts}

A canonical relation similar to $C_{CW}$ described in Thm. \ref{thm crosswell}, 
with similar geometry for $\Sigma^1,\, \Sigma^2$, and singularities types of the projections  from them,
was shown to appear in the context of synthetic aperture radar and analyzed in \cite{AFKNQ}.
The open dense subset of $C_{WA}$ described in  Thm. \ref{thm walkaway} has a similar structure.
It was shown in  \cite{AFKNQ}  that, if $A$ is an FIO of order $m$ associated with such a canonical relation,
then
\be\label{eqn comp AFKNQ}
A^*A \in I^{2m,0}(\Delta, C_1)+I^{2m,0}(\Delta, C_2)+I^{2m,0}(C_1, C_2),
\ee
where  $C_1$ is the  graph of   a canonical involution $\chi$, and $C_2$ is a two-sided fold.
It follows from \eqref{eqn ipl orders} that the order
of $A^*A$ is the same  (namely $2m$) on all three of $\Delta, C_1$ and $ C_2$,
away from their intersections,
and hence the artifacts created  by $C_1$
 and $C_2$ when attempting imaging by  backprojection are  as strong as the true image,
 and thus are nonremovable.

Due to the presumed absence of normal forms for canonical relations with this structure,
the results of  \cite{AFKNQ}, and thus
\eqref{eqn comp AFKNQ}, cannot be applied directly to the linearized scattering map
$F\in I^{-\frac12}\left(C_{CW}\right)$,
but the negative implications for artifacts are nevertheless relevant here, as can be seen by microlocalizing to 
$\Sigma^1\setminus \Sigma^2$ and $\Sigma^2\setminus \Sigma^1$,
strongly indicating but not proving the presence of strong, nonremovable artifacts in reconstructions 
from crosswell and walkaway data.
However, the presence of strong artifacts can definitely be deduced from the microlocal structure of $C_{CW}$ 
and $C_{WA}$ near points where both $\pi_L$ and $\pi_R$ are folds.
Folding canonical relations, for which both $\pi_L$ and $\pi_R$ are folds,
were first studied in the context of scattering by obstacles \cite{MeTa85}, 
and then for linearized seismics  in \cite{N00,F05,F07}.
It was shown in \cite{N00,F05} that, if $A\in I^m(C)$, then
$A^*A \in I^{2m,0}(\Delta, C_1)$  where $C_1$ is another folding canonical relation.
Since $A^*A \in I^{2m}(\Delta \setminus C_1)$ and
$A^*A \in I^{2m}(C_1 \setminus \Delta)$ by \eqref{eqn ipl orders},  the artifact created by $C_1$ is as
strong as the true image, again resulting in a nonremovable artifact.

%%%%%%%%%%%%%%%%%%%%%%%%%%%%%
%%%%%%%%%%%%%%%%% section break

\section{Appendix: Singularity classes}\label{sec sings}

\par Let $V$ and $W$ be smooth manifolds, initially of the same dimension, $n$,
and let $f : V \to W$ be a smooth function.
Let
$\Sigma(f) :=\{ x \in V : \det (df (x))=0\}$ be the set of critical points of $f$.
(This and all of the sets defined below are coordinate-independent.) The only singularities we will be concerned 
with are  those which are {\it corank one},
by which we mean points $x_0\in V$ such that
\be\label{def corankone}
\rank df(x_0)=n-1\hbox{ and } d\left(\det \left(df\right)\right)(x_0)\ne 0.
\ee
If $f$ only has corank one singularities, then $\Sigma$ is a smooth  hypersurface in $V$.
\medskip

\begin{defn}\label{def blowdown} $f:V\to W$
{\em is a} blowdown {\em if Ker} $df\subseteq T\Sigma(f)$ {\em at all points of} $\Sigma(f)$.
\end{defn}
\medskip

\begin{defn}\label{def fold}
$f$ {\em has singularities of} (Whitney) fold type
{\em if,  for
every }$x \in {\Sigma(f)}$,  ${\rm  Ker} \ df(x)$   {\em intersects } $T_x{\Sigma(f)}$  {\em transversally.}
\end{defn}

\par
\medskip

Now consider  the non-equidimensional situation.
There is some variation in the literature in terms of how these singularities are denoted.
The analogues of
Whitney folds are called {\it submersions with folds} (if $\dim(V)>\dim(W)$) or {\it cross caps} (if
$\dim(V)<\dim(W)$).
Suppose that dim $V = N$,  dim $W
= M$, with $ N \ge  M.$
For $N=M$, submersions with folds are Whitney folds,
and are denoted by $S_{1,0}$ (in the
Thom theory of $C^\infty$ singularities
\cite{GoGu})  and by $\Sigma_{1,0}$ (in the Boardman-Morin theory
\cite{mo1}) in the equidimensional case. In general,

\begin{defn}\label{def SWfold}
$f$ {\em is a} submersion  with
folds {\em if the only singularities of} $f$  {\em are of type}
$S_{1,0}$ {\em (Thom) or} $\Sigma_{N-M+1,0}$ {\em (Boardman-Morin).}
\end{defn}

\par For our purposes, we do not need to define the classes $S_{1,0}$ or $\Sigma_{N-M+1,0}$,
but simply recall that
one can verify that $f$ is a submersion with folds as follows.
At  points where
\noindent$\rank df\ge M-1$, by \cite{mo1}, we can choose suitable adapted
local coordinates on $V$ and $W$ such that $f$
has the form:
$f(x_1, x_2, \dots, x_{M-1}, x_M, \dots, x_N)=(x_1, x_2, \dots
x_{M-1}, g(x))$. The set $\Sigma(f)$ where $f$ drops rank (by $1$, by assumption) is
described by $\Sigma(f)= \{ x: \frac{\partial g}{\partial x_i}=0, \ M \leq i \leq N \}$.
Then $f$ is a submersion with folds if, for all $x\in\Sigma(f)$,
\smallskip

(i)   $\left\{ d\left(\frac{\partial g}{\partial
x_i}\right):
M \leq  i  \leq N)\right\}$  is linearly
independent (so that $\Sigma(f)$ is a
smooth submanifold of $V$); and
\smallskip

(ii) the $(N-M+1)$-dimensional kernel of
$df(x)$ is transversal to the tangent space of $\Sigma(f)$ in $T_xV$.

These conditions can be combined \cite{mo1} into
   \begin{equation}\label{swf}
\det\left[\frac{\partial^2 g }{\partial x_i
\partial x_j }\right]_{ M \leq i,j \leq N} \ne 0,
\end{equation}
and this is independent of the choice of adapted coordinates.
\medskip

\par For each $N,M$, there are a finite number of local  normal forms for a submersion
with folds, determined by the signature of the Hessian of $f$ \cite{GoGu}:

\[f(x_1,x_2, \dots, x_N) =(x_1,x_2, \dots, x_{M-1}, x_M^2 \pm
x_{M+1}^2 \pm \cdots \pm x_N^2 ).
\]
In the case relevant here, $N=4=M+1$ and the last
entry is a quadratic form in two variables.
\medskip

\par We now define the final singularity class of interest,
assuming that $f:V\to W$, with $\dim V=N < \dim W=M$.

\begin{defn}\label{def crosscap}
$f$ {\em is a} cross cap  {\em
if the only singularities of} $f$  {\em are of type} $ S_{1,0}$ (Thom) or
$\Sigma_{1,0}$ (Boardman-Morin).
\end{defn}

\par One can  identify a cross cap as follows \cite{mo1}.
At a point where  $df$ has rank $\ge N-1$, we can find suitable
adapted coordinates such that
$$f(x_1, x_2, \dots, x_{N-1}, x_N)=(x_1, x_2,
\dots, x_{N-1}, g_1, g_2, \dots g_q),$$
where $q= M-N+1$.
The set $\Sigma(f)$ where $f$ drops rank by $1$ from its maximal possible value, $N$,
is given by
$\Sigma(f) =\{x: \frac{\partial g_i}{\partial x_N}=0,  \quad  1 \leq i \leq q\}$.
Assume
that there is  an $i_0$, such that $\frac{\partial ^2
g_{i_0}}{\partial x_N^2} (0) \neq 0$. Then, $g$ has a cross cap
singularity near 0 if the map $\chi : \R^N \rightarrow \R^q $ given by
$\chi(x_1, x_2, \dots x_N)=\left(\frac{\partial g_1}{\partial x_N},
\frac{\partial g_2}{\partial x_N}, \dots, \frac{\partial
g_q}{\partial x_N}\right)$ satisfies $\rank d\chi(0) =q$. (Notice that this forces
$N \geq q$, i.e., $M \leq 2N-1$.)
These conditions can be
expressed as:
\medskip

$(i)$
$\Sigma(f)$ is smooth and of codimension $q$;
\smallskip

$(ii)$ the $N\times N$ minors of
$df$ generate the ideal of $\Sigma(f)$; and
\smallskip

$(iii)$ $\Ker (df)\cap T\Sigma(f)=(0)$.
\medskip

\par  As for folds, there is a local normal form for  cross caps, due to
\linebreak Whitney \cite{Wh45} and Morin \cite{mo1}:

\begin{equation}\label{ccnf}
f(x_1,x_2, \dots, x_N) =(x_1,x_2, \dots, x_{N-1}, x_1 x_N, \dots
x_{M-N} x_N, x_N^2).
\end{equation}

%%%%%%%%%%%%%%%%%%%%%%%%%%%%%
%%%%%%%%%%%%%%%%% section break

\section{Acknowledgements}

This paper grew out of work supported by an American Institute of Mathematics Structured Quartet Research Experience (SQuaRE).
AG was partially supported by NSF DMS-1906186.
The authors would like to thank Olga Podgornova for useful conversations.

%%%%%%%%%%%%%%%%%%%%%%%%%%%%%
%%%%%%%%%%%%%%%%% section break

\bigskip

\end{document}